\documentclass[11pt]{article}
\usepackage[letterpaper]{geometry}
\usepackage[utf8]{inputenc}

\usepackage{amsmath, amssymb, amsthm, bbm, bm}
\usepackage{authblk, enumitem, graphicx, subcaption, xcolor}
\usepackage{subfiles}
\usepackage[sc]{mathpazo}
\usepackage[linktocpage=true]{hyperref}
\hypersetup{
 colorlinks,
 linkcolor=blue,          
 citecolor=red!50!black,       
 filecolor=black,   
 urlcolor=black, 
 pdftitle={},
 pdfauthor={},
 pdfcreator={},
 pdfsubject={},
 pdfkeywords={}
}
\usepackage[sort&compress,numbers]{natbib}

\newcommand{\pro}{\mathbb{P}}

\newcommand{\R}{\mathbb{R}}
\newcommand{\N}{\mathbb{N}}

\newcommand*\diff{\mathop{}\!\mathrm{d}}
\newcommand{\pp}{\boldsymbol{p}}
\newcommand{\xx}{\boldsymbol{x}}
\newcommand{\yy}{\boldsymbol{y}}
\newcommand{\qq}{\boldsymbol{q}}

\newcommand{\cN}{\mathcal{N}}
\newcommand{\cX}{\mathcal{X}}
\newcommand{\cY}{\mathcal{Y}}
\newcommand{\cS}{\mathcal{S}}
\newcommand{\cQ}{\mathcal{Q}}

\newcommand{\keywords}[1]{\textbf{\textit{Keywords ---}} #1}

\numberwithin{equation}{section}
\newtheorem{theorem}{Theorem}[section]
\newtheorem{lemma}[theorem]{Lemma}

\newtheorem{definition}[theorem]{Definition}
\newtheorem{proposition}[theorem]{Proposition}
\newtheorem{remark}{Remark}

\begin{document}

\title{Load Balancing Under Strict Compatibility Constraints}

\author{Daan Rutten\thanks{Email: \href{mailto:drutten@gatech.edu}{drutten@gatech.edu}} }
\author{Debankur Mukherjee}
\affil{Georgia Institute of Technology}


\maketitle

\begin{abstract}

We study large-scale systems operating under the JSQ$(d)$ policy in the presence of stringent task-server compatibility constraints. Consider a system with $N$ identical single-server queues and $M(N)$ task types, where each server is able to process only a small subset of possible task types. Each arriving task selects $d\geq 2$ random servers \emph{compatible to its type}, and joins the shortest queue among them. The compatibility constraint is naturally captured by a fixed bipartite graph $G_N$ between the servers and the task types. When $G_N$ is complete bipartite, the meanfield approximation is proven to be accurate. However, such dense compatibility graphs are infeasible due to their overwhelming implementation cost and prohibitive storage capacity requirement at the servers. Our goal in this paper is to characterize the class of \emph{sparse} compatibility graphs for which the meanfield approximation remains valid.

To achieve this, first, we introduce a novel graph expansion-based notion, called \emph{proportional sparsity}, and establish that systems with proportionally sparse compatibility graphs match the performance of a fully flexible system, asymptotically in the large-system limit. Furthermore, for any $c(N)$ satisfying 
$$\frac{Nc(N)}{M(N)\ln(N)}\to \infty\quad \text{and}\quad c(N)\to \infty,$$
as $N\to\infty$, we show that proportionally sparse random compatibility graphs can be designed, so that the degree of each server is at most $c(N)$. This reduces the server-degree almost by a factor $N/\ln(N)$, compared to the complete bipartite compatibility graph, while maintaining the same asymptotic performance. Extensive simulation experiments are conducted to corroborate the theoretical results. 
\end{abstract}



\section{Introduction}

\subsection{Background and motivation}

A canonical model for large-scale systems, such as data centers and cloud networks, consists of a large number of parallel servers with dedicated queues. 
Tasks arrive into the system sequentially in time and are immediately and irrevocably assigned, using some efficient load balancing algorithm, to one of these queues, where they wait until executed.
Due to ever-increasing heterogeneity in the incoming traffic, these systems typically suffer from stringent task-server compatibility constraints. 
Indeed, executing a task at a server requires some pre-stored data, and being able to serve all possible task types comes with an excessive storage capacity requirement~\cite{WZYTZ16, XYL16} and an overwhelming implementation complexity~\cite{MHCD10, RTGKK12, TX17}. 
Consequently, full flexibility in task allocation is not a luxury large-scale systems can afford.
It is therefore important to understand the performance of load balancing algorithms under sparser compatibility constraints, where tasks of a particular type can only be served by a relatively small number of servers, naturally viewed as neighbors in a bipartite compatibility graph between the servers and the task types.

The analysis of load balancing algorithms for large-scale systems dates back to the seminal works
by Vvedenskaya et al.~\cite{VDK96} and Mitzenmacher~\cite{Mitzenmacher96}.
Since then, using meanfield techniques, there has been significant progress in our understanding of the performance of various algorithms. 
However, many of these heuristics turn out to be false in the presence of compatibility constraints.
A widely studied algorithm in this area is the JSQ($d$) or the `power-of-$d$' scheme, where each arriving task is assigned to the shortest of $d$ randomly selected queues. 
The JSQ($d$) scheme is popular for its low-complexity implementation and excellent delay performance.
However, an ill-designed compatibility structure can even lead to instability or poor delay performance of a system operating under the JSQ($d$) scheme.
In spirit, the nature of this observation is similar to the famous Braess's paradox~\cite{Braess68, RT02} in networks.
We discuss this example in more detail in Remark~\ref{rem:nonmonotone} below.

The lack of a thorough understanding of large-scale systems with compatibility constraints may be attributed to the scarcity of the theoretical toolbox to analyze such systems.
Performance analysis of large-scale systems has flourished in the last three decades due to the abundance of sophisticated meanfield techniques and in particular, the asymptotic analysis of density-dependent population processes~\cite{EK2009}.
This has provided a firm theoretical basis to analyze these systems.
In the presence of arbitrary compatibility constraints, the servers become non-exchangeable, which breaks a core assumption that lies at the foundation of the classical meanfield framework. 
Our goal in this paper is to understand the effect of compatibility constraints on the performance of large-scale systems and in particular, characterize a large class of sparse compatibility graphs that match the performance of a fully flexible system asymptotically in the large-system limit.
In doing so, we will also make progress in developing new approaches to analyze such structurally constrained large-scale systems driven by stochastic inputs.

\subsection{Our contributions}
Consider a system with $N$ single-server queues and $M(N)$ task types. 
We will be looking at a scaling regime where both $N, M(N)\to \infty$.
The task-server compatibility is captured in terms of a bipartite graph $G_N$ between the servers and task types.
That is, a server $i$ shares an edge in $G_N$ with a task type~$j$, if task type~$j$ can be processed by server $i$.
Following the JSQ($d$) policy, when a task of type $j$ arrives, $d$ servers which share an edge with $j$ are sampled uniformly at random and the task is routed to the shortest of the sampled queues.
The quantity of interest is the global occupancy process $\qq^N(t) = (q_1^N(t), q_2^N(t)\ldots)$, where $q_i^N(t)$ denotes the fraction of servers with queue length at least $i$ at time $t$ in the $N$-th system.
Note that the case when $G_N$ is complete bipartite corresponds to the fully flexible system.
Our focus is to identify the sparsest compatibility structures that preserve the performance benefits of a fully flexible system, asymptotically as $N\to \infty$.
In other words, we study the \emph{sparsity condition for the compatibility graph, which preserves the validity of the meanfield approximation}.

It is reasonable to guess that well-designed compatibility graphs with `sufficient' amount of expansion property should preserve the effects of full flexibility, asymptotically, in the large-system limit. 
However, 
identifying the right notion of expansion and thereby establishing precise limit laws for process-level and steady state occupancy processes remains a notoriously challenging problem, and has inspired several research works, as we will discuss below. 
In this work we attempt to make progress in this direction by developing 
new approaches to tackle the non-exchangeability.
Specifically, our results can be categorized into two groups.\\

\noindent
\textbf{(1) Arbitrary deterministic compatibility graphs.} 
We start by considering an arbitrary deterministic sequence of graphs $\{G_N\}_{N\geq 1}$, indexed by the number of servers $N$, and define a novel notion of expansion, which we call \emph{proportional sparsity}, see Definition~\ref{con:graphseq}. We show that if the sequence of compatibility graphs is proportionally sparse, then as $N\to\infty$, on any finite time interval, the occupancy process $\qq^N(\cdot)$ under the JSQ($d$) policy converges to the same meanfield limit as the sequence of fully flexible systems. 
In fact, this process-level limit result extends to a broad class of load balancing algorithms, for which the assignment decision depends `smoothly' on the empirical queue length distribution of the compatible servers. We call such algorithms \emph{Lipschitz continuous task assignment policies}, see Definition~\ref{con:lip-routing}.
An important step to prove the process-level limit is to show that for almost all dispatchers, the empirical queue length distribution observed in its neighborhood, is close to the empirical queue length distribution observed among all servers in the system. This allows us to construct a coupling between the constrained system and the fully flexible system and establish that the $\ell_1$-distance between the global occupancy processes in two systems is small uniformly over any finite time interval.

For the interchange of limits and hence the convergence of steady state, two more key ingredients that we need are ergodicity of the prelimit system (for each fixed $N$) and the tightness of steady states in an appropriate sense. 
Note that if $G_N$ is not complete bipartite, the occupancy process $\qq^N(\cdot)$ is no longer Markovian.
Consequently, one needs to be careful in defining its time asymptotics and hence, the interchange of limits.
For ergodicity of the underlying Markov process, we need the graph sequence to satisfy a certain \emph{subcriticality condition} that was first introduced in~\cite{Bramson11}, see Definition~\ref{def:subcritical}. The tightness, however, is technically more challenging. In particular, we need to show that the sequence of steady state occupancy is tight with respect to a certain weighted $\ell_1$ norm. For this, we construct a collection of Lyapunov functions, which provide uniform tail bounds on the steady state of the global occupancy process.

Combining the above results, we conclude in Theorem~\ref{th:maintheorem} that if a sequence of graphs is proportionally sparse and satisfies the subcriticality condition, then both finite-time dynamics and steady state behavior of the empirical queue length process coincide with that of a fully flexible system, asymptotically as $N\to\infty$.
It is worth highlighting that in the above interchange of limits, we do not impose any restrictions on how the number of task types $M(N)$ scales with $N$. This includes the two popular scenarios $M(N)=$ constant and $M(N)=N$ as special cases.\\

\noindent
\textbf{(2) Random compatibility graphs.}
The results for deterministic graph sequence provide us all the theoretical framework needed to analyze these systems.
Next, we exploit these results to construct random sparse compatibility graphs with desired performance benefits.
In the context of data-file placement or content replication in large-scale systems, the degree of a server in the compatibility graph can be thought to be roughly proportional to the storage capacity requirement of that server. It is also considered to be a measure of complexity of the network. To this end, we consider two cases.

First, suppose that the servers are constrained to have degrees exactly equal to $c(N)$. In this case construct $G_N$ by selecting $c(N)$ task types for each server, independently uniformly at random, without replacement, from the set of all task types. For such a randomly constructed compatibility graph, we establish that the empirical queue length distribution of the system has the same asymptotic law as the fully flexible system, both process-level and in steady state, if $c(N)\gg M(N)\ln(N)/N$ and $c(N)\gg 1$, see Theorem~\ref{th:drg} for details.

Second, we consider a system that allows for inhomogeneous levels of flexibility for different task types. In this case, the compatibility graph is constructed by selecting each edge incident to a task-type $w \in W_N$ with probability $p_w(N)$, independently of other edges. Thus, task type $w$ will have an average degree $Np_w(N)$. In this case, we show that the empirical queue length distribution of the system has the same asymptotic law as the fully flexible system, both process-level and in steady state, if $\min_{w \in W_N} p_w(N)$ and the $\ell_2$ norm of the inverse probability vector $(1/p_w(N))_{w\in W_N}$ satisfy suitable growth conditions, see Theorem~\ref{th:erg} for details.

To prove the results for random instances, we verify using concentration of measure arguments, that the graph sequence satisfies both the proportional sparsity and the subcriticality conditions, under the respective growth rate conditions as $N\to\infty$.

\subsection{Related works}
The effect of flexibility in the task assignment in large-scale systems was first studied by Turner~\cite{Turner98}, who considered two types of arriving customers, those that have no routing choice, and those that employ the JSQ($d$) strategy. It was shown that even a small amount of routing choice can lead to substantial gains in performance through resource pooling. 

Relatively recently, there have been a number of works analyzing load balancing algorithms for large-scale systems, where queues themselves are interconnected by some graph topology. In these models, each queue has an independent dedicated stream of external arrivals and each external arrival must be assigned instantaneously and irrevocably to one of the neighboring queues including the one where it first appeared. While these models cannot directly be used to capture the task-server compatibility constraints, from a high level, they can be viewed as an undirected graph version of our model in the special case when $M(N)=N$, see \cite[Remark 4]{BMW17}.
In this line of works, motivated by the bike-sharing network, Gast~\cite{Gast15} studies a system of queues connected by a ring topology. To deal with the long-range dependencies among the queue-length processes arising from the restricted graph topology, the work proposes a pair-approximation to describe the steady state system. 
When the ordinary JSQ policy is used at each vertex, that is, when each arriving task joins the shortest of all the neighboring queues, Mukherjee et al.~\cite{MBL17} develop a coupling-based approach to establish criteria for asymptotic optimality on fluid and diffusion scale.
A key ingredient in this approach is the monotonicity of the system with respect to edge addition. 
That is, performance of a system gets better, in the sense of stochastic majorization of the occupancy process, if more edges are added to the underlying graph.
The scenario becomes fundamentally more challenging, however, if the system lacks the above-mentioned monotonicity with respect to edge addition. 
One such scenario is when the JSQ($d$) policy is considered at each vertex, instead of the JSQ policy.
That is, each task is assigned to the shortest queue among the one it first appears and its $d-1$ randomly selected neighbors.
Even a first-order property such as stability is non-trivial in this case. See Remark~\ref{rem:nonmonotone} for a related illustration of this non-monotonicity. 
It is this non-monotonicity that makes the scenario considered in the current article very different from those in the state-of-the-art literature.
In fact, this non-monotonicity hints that \emph{expansion properties that are monotone with respect to edge addition cannot provide sufficient criteria for getting the same asymptotic limit law as the fully flexible system}.

Contemporaneously to the current article, Weng et al.~\cite{WZS20} consider the join-the-fastest-shortest-queue (JFSQ) and the join-the-fastest-idle-queue (JFIQ) policy for systems with task-server compatibility constraints, where the arrival rates of the task types and the service rates of the servers are heterogeneous. 
Specifically,~\cite{WZS20} obtains finite-system bounds on the mean response time, and generalizing the Lyapunov drift
method, shows that under a `well-connected' graph condition, the JFSQ and JFIQ policies can achieve the minimum steady-state response time in both the many-server regime and the sub-Halfin-Whitt regime (when the system load approaches one at a suitable rate), asymptotically as~$N\to\infty$.

In a work by Budhiraja et al.~\cite{BMW17}, sufficient conditions on the graph sequence are obtained, so that the queue length process under the JSQ($d$) policy has the same fluid limit on any finite time interval as the complete graph.
Their method relies on an asymptotic coupling of the queue length process with an infinite-dimensional McKean-Vlasov process, which does not easily generalize to steady state. Indeed, as mentioned in \cite[Section 4]{BMW17}, even the existence of a time asymptotic limit was not clear for this system.

In another interesting line of research, Tang and Subramanian~\cite{TS19a, TS19}  analyze a variant of the classical JSQ($d$) policy, where the $d$ servers are sampled through $d$ independent non-backtracking random walks on a high-girth graph.
The motivation here is to reduce the amount of randomness used in implementing the classical JSQ($d$) policy. \\

There has been a rich literature in the stability analysis of load balancing algorithms of finite-sized systems. 
Related to the output-queued model considered in this paper, the stability analysis dates back to Stolyar~\cite{Stolyar95, Stolyar05} and Chernova and Foss~\cite{FC98}.
Building on the framework of~\cite{FC98}, Bramson~\cite{Bramson11} analyzed stability of JSQ-type systems under a broad class of policies, including the JSQ($d$) policy, where the service times and inter-arrival times follow general distributions. Instead of a bipartite compatibility graph, in this work, there is an independent arrival stream of tasks of rate $\lambda_S$ corresponding to each subset $S \subseteq [N]$ of servers. Tasks in the arrival stream $S$ join the shortest queue among $S$. Bramson~\cite{Bramson11} proposes a sufficient subcriticality condition on the arrival rates $\lambda_S$ and shows that the system is ergodic under this condition.
We will use the above subcriticality condition to establish ergodicity of the system for each fixed $N$.
However, the results of~\cite{Bramson11} do not guarantee that the steady-state workload in the system scales as $\Theta(N)$ as $N\to \infty$.
This is crucial in the large-system limit, since it relates to the tightness of the sequence of steady state occupancy.
For this, we use the Lyapunov function approach, as in~\cite{WMSY18, WMSY18a}, and establish moment bounds~\cite{Hajek82, MT93} to obtain uniform bounds on the tail of the stationary occupancy process.

More recently, Cardinaels et al.~\cite{CBL19} analyze stability conditions and obtain performance bounds of a general model for load balancing with affinity relations.
In this setup, each arriving task can be routed to either a fast, primary selection of servers or a secondary selection with a slower processing speed. 
Cruise et al.~\cite{CJS20} establish stability for a similar problem where the task-server constraints are modeled as a hypergraph. 
In the area of redundancy scheduling under compatibility constraints, Cardinaels et al.~\cite{CBL20} study the case when each task may only be replicated to a specific set of servers described by a compatibility graph. In the classical heavy-traffic regime (fixed number of servers and load approaches the boundary of the capacity region) and under appropriate conditions on the graph,~\cite{CBL20} establishes that the system with graph-based redundancy scheduling operates as a multi-class single-server system.\\

On the scheduling side, there has been significant development in the analysis of multiserver input-queued systems with multiple task types; see~\cite{TX13, TX17, GS12, HAP19, yekkehkhany2018gb, yekkehkhany2020blind, AHP12} and the references therein. 
In this area, the work that is closest, in spirit, to our setup, is by Tsitsiklis and Xu~\cite{TX13, TX17}. Here, an input-queued system is considered, where $N$ servers are connected to $r N$ queues by a bipartite compatibility graph, where $r>0$ is a fixed constant that does not depend on~$N$. Each queue~$i$ receives an independent arrival stream of rate $\lambda_i$ and tasks remain in the queue until a server becomes available. In this setup, Tsitsiklis and Xu~\cite{TX13, TX17} establish that if the average degree of the queues $c(N) \gg \ln(N)$, then there exists a family of expander-graph-based  flexibility architectures  and a scheduling policy that stabilizes almost all admissible arrival rates and is throughput optimal. \\

Lastly, this work also fits into the recent line of works on load balancing for systems with multiple dispatchers~\cite{BBL17a, Stolyar17, ZSW20, VKO20}. The analysis in these cases is often more challenging than the classical setup. However, strict task-server compatibility are typically not considered in these works.
We refer to~\cite{BBLM18} for a recent survey on the load balancing algorithms.

\subsection{Notation and organization}
The remainder of the paper is organized as follows.
In Section~\ref{sec:model} we describe the model in detail and introduce notations related to the underlying Markov chain and its state space.
Section~\ref{sec:mainresults} lists the main theorems and discusses their ramifications. 
Section \ref{sec:proofoverview} provides the proofs of the main theorem involving the deterministic graph sequence.
Due to space restriction, proofs of statements marked $(\bigstar)$ have been omitted; they can be found in the appendix.
In Section~\ref{sec:simulation} we present simulation experiments, both to support the analytical results and to examine the performance of systems with compatibility constraints that are not analytically tractable.  
Finally, Section~\ref{sec:conclusion} summarizes our results and discusses directions for future research.\\

A complete bipartite graph with $N$ servers and $M(N)$ dispatchers will be denoted by $K_{N,M}$. We denote by $\ell_1$, the normed vector space with norm $\lVert \qq \rVert_1 = \sum_{i = 1}^\infty \lvert q_i \rvert$.
For some positive sequence $\boldsymbol{\omega} = (\omega_1, \omega_2,\ldots)$, denote by $\ell_1^\omega$, the normed vector space with the weighted-$\ell_1$ norm $\lVert \qq \rVert_1^\omega = \sum_{i = 1}^\infty \omega_i \lvert q_i \rvert$.
For any set $V$, $|V|$ denotes its cardinality. 
We adopt the usual notations $O(\cdot), o(\cdot), \omega(\cdot),$ $\Omega(\cdot)$, and $\Theta(\cdot)$ to describe asymptotic comparisons.
For two positive deterministic sequences $(f(n))_{n\geq 1}$ and $(g(n))_{n\geq 1}$, we write $f(n)\ll g(n)$ (respectively, $f(n)\gg g(n)$), if $f(n) = o(g(n))$ (respectively, $f(n) = \omega(g(n))$).

\section{Model description}\label{sec:model}
\begin{figure}
  \centering
  \includegraphics[width=0.4\textwidth]{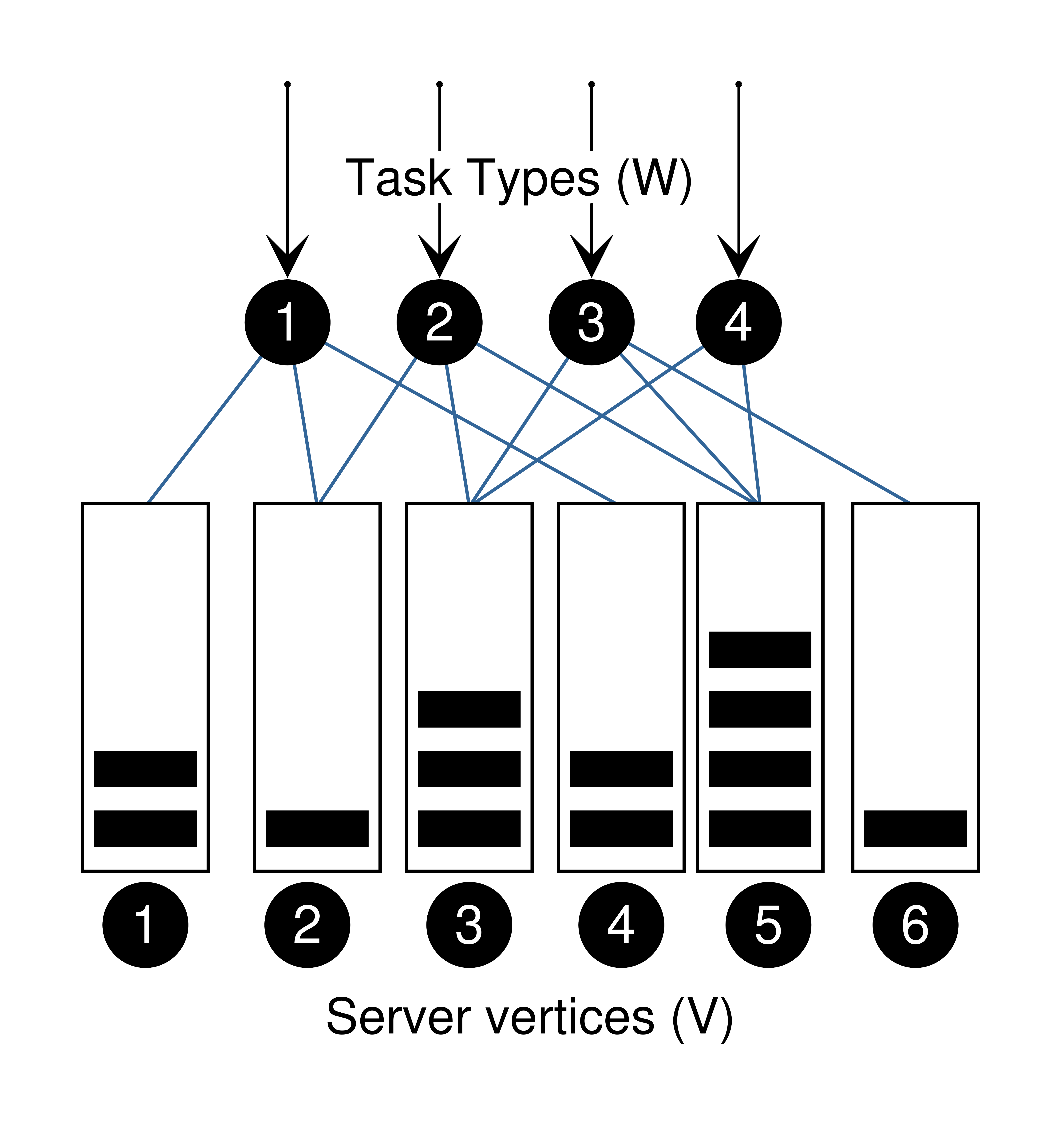}
  \caption{A schematic overview of the system with task types $W_N$, servers $V_N$, and their compatibility relation.}
  \label{fig:dispatchers}
\end{figure}

We consider a system that consists of a set of dispatchers $W_N$ and a set of servers $V_N$, with $\lvert V_N \rvert = N$ and $\lvert W_N \rvert = M(N)$. 
Each dispatcher handles arrivals of a particular task type. 
Task-type $w \in W_N$ can only be served by a subset $\cN_w \subseteq V_N$ of servers.
The set of servers $\cN_w$ compatible to the task type $w$, can naturally be viewed as neighbors of $w$ in
a bipartite graph $G_N = (V_N, W_N, E_N)$ between $V_N$ and $W_N$, where $E_N \subseteq V_N \times W_N$ is the set of edges.
In other words, $w\in W_N$ and $v\in V_N$ share an edge in $G_N$, if server $v$ has the resources required to process task-type $w$; see Figure~\ref{fig:dispatchers} for an illustration of the model.
Without loss of generality, we assume $G_N$ to be connected.
Note that the fully flexible system corresponds to $G_N$ being complete bipartite.
Each server has a dedicated queue with infinite buffer capacity and operates under a non-idling service discipline that is oblivious to the actual service requirements (e.g., FCFS). Below we will interchangeably use the terms \emph{task-type vertex} and \emph{dispatcher}.

Each dispatcher receives an external stream of arrivals as a Poisson process of rate $\lambda N / M(N)$, independently of the other processes.
Thus, the total arrival rate of the system is $\lambda N$. 
The service times of the tasks are exponentially distributed with unit mean, independently of each other.
For the stability of the system, assume $\lambda < 1$. 
Each dispatcher employs the JSQ($d$) policy. That is, when a task arrives at $w \in W_N$, the dispatcher samples queue lengths at $d\geq 2$ servers uniformly at random from its neighborhood $\cN_w$. The task is assigned to the shortest of the sampled queues.\\


At time $t\geq 0$, let $I^N_v(t)$ denote the queue length of server $v\in V_N$. Also, let $\xi_{vw}^N$ denote the edge occupancy in $G_N$, that is, for $v\in V_N$ and $w\in W_N$, $\xi_{vw}^N=1$ if $(v,w)\in E_N$ and $\xi_{vw}^N=0$ otherwise.
Note that $\big(I_v^N(\cdot), \xi_{vw}: v\in V_N, w\in W_N\big)$ is a Markovian state descriptor of the system.
Denote this Markov process at time $t\geq 0$ by $\Phi(G_N, t)$ and the state space by $\cS_N$.
Since for each system, we keep the graph structure fixed beforehand, we leave the edge occupancy implicit in the state and think of $\Phi(G_N, t)$ as the vector of tagged queue lengths $(I_v^N(\cdot): v\in V_N)$ and $\cS_N$ as~$\mathbb{N}^N$.


We introduce a few shorthand notations. For a state $z \in \cS_N$, let $\cX_i(z) \subseteq V_N$ denote the set of servers with queue lengths exactly equal to $i$, $X_i(z) := \lvert \cX_i(z) \rvert$ and $x_i(z) := X_i(z) / N$. 
Also, let $\xx(z) := (x_i(z))_{i \geq 0}$. We refer to $\xx(z)$ as the \textit{global empirical queue length distribution} (GEQD). 
Furthermore, let $X_i^w(z) := \lvert \cX_i(z) \cap \cN_w \rvert$ denote the number of servers with queue length exactly equal to $i$  in the neighborhood $\cN_w$ of $w \in W_N$, and  $x_i^w(z) := X_i^w(z) / \lvert \cN_w \rvert$ and $\xx^w(z) := (x_i^w(z))_{i \geq 0}$. 
We refer to $\xx^w(z)$ as the \textit{local empirical queue length distribution} (LEQD).
The space of (global and local) empirical queue length distributions is denoted as 
$$\cX := \Big\{ \xx \in \ell_1\ \Big|\ \sum_{i = 0}^\infty x_i = 1 \text{ and } x_i \geq 0 \; \forall i \in \mathbb{N} \Big\}.$$

For a state $z \in \cS_N$, let $\cQ_i(z)$ denote the set of servers with queue length \textit{at least} $i$, $Q_i(z) := \lvert \cQ_i(z) \rvert$ and $q_i(z) := Q_i(z) / N$ such that $q_i(z) = \sum_{j = i}^\infty x_j(z)$. Also, let $\qq(z) := (q_i(z))_{i \geq 1}$. 
We refer to $\qq(z)$ as the \textit{global occupancy}.
Also, for $w \in W_N$, let $Q_i^w(z) := \lvert \cQ_i(z) \cap \cN_w \rvert$ denote the number of servers with queue length at least $i$ in $\cN_w$, $q_i^w(z) := Q_i^w(z) / \lvert \cN_w \rvert$, and $\qq^w(z) := (q_i^w(z))_{i \geq 1}$. We refer to $\qq^w(z)$ as the \textit{local occupancy}.
Local and global occupancy take values in the space
$$\cY := \Big\{ \qq \in \ell_1\ \Big|\ q_i\in [0,1], q_i\geq q_j, i<j, i,j\in \mathbb{N} \Big\}.$$

\section{Main results}
\label{sec:mainresults}
\subsection{Arbitrary deterministic graphs}

We start by introducing the notion of proportional sparsity and the subcriticality condition for a deterministic sequence of bipartite graphs $\{ G_N \}_{N \geq 1}$ and discuss their ramifications. Proportional sparsity provides a sufficient expansion property of the compatibility graph so that, on any finite time interval, the occupancy process of systems under a broad class of task assignment policies has the same weak limit as the fully flexible system. The subcriticality condition bounds the maximum load on any server and  implies that the underlying Markov process is ergodic and the steady state global occupancy is tight in the appropriate space. Together, the proportional sparsity and the subcriticality condition imply the interchange of limits for the global occupancy process.\\

\noindent
\textbf{Proportional sparsity.}
The condition of proportional sparsity requires the edges in the bipartite graph to be fairly distributed in an appropriate sense.
\begin{definition}
\label{con:graphseq}\textit{
(Proportionally sparse graph sequences)
Let $G_N = (V_N, W_N, E_N)$ be a sequence of connected graphs indexed by the number of servers $|V_N| =N$. 
The sequence $\{G_N\}_{N\geq 1}$ is called proportionally sparse if for each $\varepsilon > 0$,
\begin{equation}
    \sup_{U \subseteq V_N} \left\lvert \left\{ w \in W_N \mid \left\lvert \frac{\lvert \cN_w \cap U \rvert}{\lvert \cN_w \rvert} - \frac{\lvert U \rvert}{N} \right\rvert \geq \varepsilon \right\} \right\rvert / M(N) \to 0, \quad\text{as}\quad N \to \infty.
\end{equation}
}
\end{definition}

\begin{remark}\normalfont
From a high level, the class of proportionally sparse graph sequences contains all graphs obtained after \emph{two-step sparsification} of the complete bipartite graph.
To see this, note that the complete bipartite graph is the only graph for which $\frac{\lvert \cN_w \cap U \rvert}{\lvert \cN_w \rvert} = \frac{\lvert U \rvert}{N} $ for all $U\subseteq V_N$ and $w\in W_N$.
Now, the first step of sparsification allows for a wiggle-room of $\varepsilon$, however small, in the above difference for all $U$, and the second step allows for $o(N)$-many dispatchers to have the above difference larger than $\varepsilon$.
As we will see, after these two steps of sparsification, the class of graph sequences that satisfy this property will be large, and will contain graph sequences that are much sparser than the complete bipartite graph.
\end{remark}


\begin{remark}[Proportional sparsity and quasi-randomness]
\normalfont\label{rem:quasirandom}
In the dense case when $|\cN_w|$ is $\Theta(N)$, the definition of proportional sparsity is related to the notion of quasi-random bipartite graphs. To obtain a random server network one approach is to construct a random graph and use its structure in the network. An alternative approach is to take a deterministic graph and question whether it is sufficiently `random'. For a sequence of these deterministic graphs it is possible to verify whether it satisfies certain properties that random graphs are expected to have. A sequence of graphs satisfying such properties is called \textit{quasi-random}.
This notion was proposed in a seminal paper by Chung, Graham, and Wilson~\cite{CGW89}, and has subsequently been used in developing numerous algorithmic heuristics. 
Lemma \ref{lem:quasi-to-prop} below states that quasi-randomness implies proportional sparsity.
\begin{lemma}[$\bigstar$]\label{lem:quasi-to-prop}
If $\{G_N\}_{N\geq 1}$ is a quasi-random sequence of graphs, it must be a proportionally sparse graph sequence.
\end{lemma}
\noindent
The proof of Lemma \ref{lem:quasi-to-prop} is provided in Appendix \ref{app:proofquasirandom}. However, even in the dense case, the \emph{converse of Lemma~\ref{lem:quasi-to-prop} is not true}. 
The main reason is the inherent symmetry assumption of quasi-randomness, whereas a proportionally sparse graph sequence can have very inhomogeneous degrees.
To see this, we refer to Theorem~\ref{th:erg}, which states that a broad class of sequences of inhomogeneous random graphs are proportionally sparse. 
\end{remark}

\noindent
\textbf{Subcriticality condition.}
We continue by introducing a condition on the maximum load at any server. The subcriticality property will allow us to prove ergodicity of the system and tightness of the steady state occupancy process.
\begin{definition}\label{def:subcritical}\textit{
(Subcriticality condition)
Let $G_N = (V_N, W_N, E_N)$ be a graph sequence. The sequence $\{G_N\}_{N\geq 1}$ is said to satisfy the subcriticality condition if for all $N \geq 1$, $w \in W_N$ and $U \subseteq V_N$ there exists a probability distribution $\gamma_w^U(\cdot)$ on $U$ such that
\begin{equation}
\label{eq:subcritical}
    \limsup_{N \to \infty} \max_{v \in V_N} \frac{N}{M(N)} \sum_{w \in W_N} \binom{\lvert \cN_w \rvert}{d}^{-1} \sum_{\substack{U \subseteq \cN_w \\ \lvert U \rvert = \min(d, \lvert \cN_w \rvert)}} \gamma_w^U(v) \leq 1.
\end{equation}
}
\end{definition}

The probability distribution $\gamma_w^U(v)$ can be interpreted as a \emph{static} randomized task assignment policy for tasks arriving at dispatcher $w \in W_N$, when the subset $U \subseteq V_N$ is selected as the $d$ chosen servers. To this end, to understand the subcriticality condition intuitively, think of a new system where the task allocation is done as follows:  
when a task arrives at a dispatcher $w \in W_N$, similar to the JSQ($d$) policy, a subset $U$ of size $d$ is sampled from the neighborhood $\cN_w$ of $w$ (if the neighborhood $\cN_w$ contains less than $d$ servers, then set $U = \cN_w$). The task is then routed to a server $v \in U$ with probability $\gamma_w^U(v)$.
Since each dispatcher receives a load of $\lambda N / M(N)$, the stability condition requires that $\lambda^{-1}$ times the load received by any server under this static task assignment policy $\gamma_w^U(v)$ to be less than one. It is thus sufficient for stability to find  one such static task assignment policy, depending on the underlying graph, satisfying this condition. 

\begin{remark}\normalfont
As mentioned in the introduction, Bramson~\cite{Bramson11} and more recently, Cardinaels et al.~\cite{CBL19} analyzed the stability of JSQ-type policies in a related framework. The subcriticality condition~\eqref{eq:subcritical} is equivalent to the those stated in \cite{Bramson11, CBL19}, although our condition is stated for a sequence of graphs, later to be used for tightness of the steady state occupancy process, and the latter conditions are stated for a fixed system. As noted in~\cite[page 1571]{Bramson11}, the subcriticality condition as stated in Definition~\ref{def:subcritical}, is \emph{also a necessary condition} for stability.
\end{remark}

\begin{remark}[Non-monotonicity]
\label{rem:nonmonotone}
\normalfont
The satisfiability of neither the proportional sparsity property nor the subcriticality condition, are monotone with respect to the addition of edges in the compatibility graph. 
As an example, consider the graph in which each dispatcher is perfectly matched to exactly one server, as in Figure \ref{fig:nonmonotone1}. This graph satisfies the subcriticality condition. Now alternatively, consider the graph in Figure \ref{fig:nonmonotone2}. The latter graph contains all the edges of the graph in Figure \ref{fig:nonmonotone1}, yet, it is not hard to verify that it does not satisfy the subcriticality condition for $d = 2$.
This is because for any task arriving at dispatchers 3, 4, 5, and 6, it will be assigned to either server 1 or 2 with probability $1/3$.
This makes the load on servers 1 and 2 higher than 1, under any static task assignment policy. Similarly, the notion of proportional sparsity is also not monotone in the addition of edges. For example, adding a disproportionate number of edges between only a subset $U \subset V_N$ and the dispatchers $W_N$ will invalidate proportional sparsity with respect to that set $U$. 
\begin{figure}
\centerline{\begin{subfigure}[t]{.33\textwidth}
    \centering
    \includegraphics[width=\textwidth]{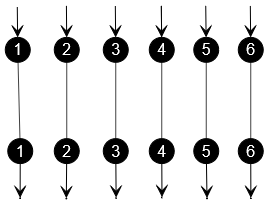}
    \caption{A simple compatibility graph which satisfies the subcriticality condition.}
    \label{fig:nonmonotone1}
\end{subfigure}\hspace{5mm}%
\begin{subfigure}[t]{.33\textwidth}
    \centering
    \includegraphics[width=\textwidth]{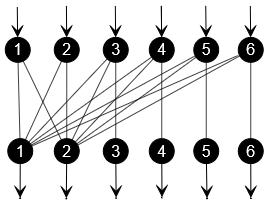}
    \caption{A compatibility graph with increased flexibility, which does \emph{not} satisfy the subcriticality condition for~$d = 2$.}
    \label{fig:nonmonotone2}
\end{subfigure}}
\caption{An example illustrating that increasing flexibility may not always lead to better performance.}
\end{figure}
\end{remark}


Now we have all the ingredients to state the main result for deterministic sequences of compatibility graphs.

\begin{theorem}
\label{th:maintheorem}
Let $\{ G_N \}_{N \geq 1}$ be a proportionally sparse graph sequence.
Then, on any finite time interval $[0, T]$, the scaled occupancy process $\qq(\Phi(G_N, t)) = (q_1(\Phi(G_N, t)), q_2(\Phi(G_N, t)), \dots)$ converges weakly with respect to the Skorohod-$J_1$ topology, as $N\to\infty$, to the process $\qq^*(t) = (q_1^*(t), q_2^*(t), \dots)$, given by the unique solution of the system of ODEs:
\begin{equation}\label{eq:ODE-gen}
    \frac{d q_i^*(t)}{dt} = \lambda (q_{i - 1}^*(t)^d - q_i^*(t)^d) - (q_i^*(t) - q_{i + 1}^*(t)), \quad\text{for } i = 1, 2, \dots,
\end{equation}
provided $\qq^*(0) \in \cY$ and $\lVert \qq(\Phi(G_N, 0)) - \qq^*(0) \rVert_1 \to 0$ as $N \to \infty$.

Moreover, if the graph sequence $\{G_N\}_{N \geq 1}$ satisfies the subcriticality condition, then for each fixed $N$, the Markov process is ergodic and $\qq(\Phi(G_N, \infty)) $ converges weakly to $\qq^*(\infty)$ as $N \to \infty$, where $\Phi(G_N, \infty)$ is a random variable distributed as the steady state of $(\Phi(G_N, t))_{t\geq 0}$ and
\begin{equation}\label{eq:qinf}
    q_i^*(\infty) = \lambda^{\frac{d^i - d}{d - 1}}, \quad\text{for }i = 1, 2, \dots.
\end{equation}
\end{theorem}

The proof of Theorem \ref{th:maintheorem} is provided in Section \ref{sec:proof-main}.


\begin{remark}\normalfont
As mentioned before, the process-level convergence result in Theorem~\ref{th:maintheorem} holds for a much broader class of assignment policies, when the assignment decision depends \emph{smoothly} on the local empirical queue length distribution (LEQD). We call this class the \emph{Lipschitz continuous} task assignment policies, which is introduced in Definition \ref{con:lip-routing} below in Section \ref{sec:overviewtransient}.
\end{remark}

\begin{remark}\normalfont
The system of ODEs in~\eqref{eq:ODE-gen} can be recognized as the meanfield limit of the classical JSQ($d$) policy~\cite{Mitzenmacher96, VDK96}.
Thus, Theorem~\ref{th:maintheorem} extends the validity of meanfield approximation for the class of proportionally sparse graph sequences that satisfy the subcriticality condition.
In other words, as we will see in the next section, Theorem~\ref{th:maintheorem} shows that the performance benefits of the fully flexible system can be preserved while making the compatibility graph significantly sparser.
\end{remark}


\subsection{Randomly designed compatibility graphs}

We have seen two sufficient conditions on deterministic graph sequences to establish asymptotically equivalence of the JSQ($d$) policy under limited and full flexibility in task allocation. Given a graph sequence, the conditions can be verified. 
This section will provide two simple ways of constructing a random compatibility graph, both satisfying the two conditions almost surely, in the large-system limit.
Note that although the graph is random in the two cases below, once constructed, it remains fixed for the system.
That is, in the terminology of random processes in random environment, we obtain a quenched limit theorem.
Below, the almost sure statements involving the sequence of random graphs $\{G_N\}_{N\geq 1}$, are with respect to any probability
measure $\mathbb{P}_0$ on $\prod_{N} \{0,1\}^{N \times M(N)}$, such that its
projection on $\{0,1\}^{N \times M(N)}$ corresponds to the distribution of the adjacency matrix of $G_N$.\\

\noindent
\textbf{Hard constraint on server-degrees.}
As mentioned before, the degree of the servers in the compatibility graph is an important measure of sparsity, as it is roughly proportional to the required storage capacity.
For that reason, we will first consider the case when all the servers have degree exactly equal to some fixed number $c(N)$, that is much smaller than $M(N)$, the server-degree for a fully flexible system.
\begin{theorem}[$\bigstar$]
\label{th:drg}
Let $c(N)\leq M(N)$ be a sequence of positive integers satisfying 
$$c(N) \to \infty\quad\text{and}\quad \frac{Nc(N)}{M(N)\ln(N)}\to \infty,\quad \text{as}\quad N\to\infty.$$ 
Also, construct $G_N$ as follows: for each $v\in V_N$, select $c(N)$ edges from $\{ (v, w) \in V_N \times W_N \mid v \in V_N \}$ uniformly at random, without replacement.
Then the sequence $\{G_N\}_{N\geq 1}$ is proportionally sparse and satisfies the subcriticality condition, almost surely.
Consequently, the conclusions of Theorem~\ref{th:maintheorem} hold.
\end{theorem}
The proof of Theorem \ref{th:drg} is technically involved. It uses concentration of measure arguments repeatedly to establish structural properties of the compatibility graphs, and is provided in Appendix~\ref{app:drg}.

\begin{remark}\normalfont
Observe that Theorem~\ref{th:drg} guarantees that the validity of the meanfield approximation can be retained asymptotically, by uniformly reducing the server-degrees by almost a factor of $N / \ln(N)$, compared to the fully flexible system where the degree of each server is $M(N)$.
Also, note that if $M(N) = O(N / \ln(N))$, then the convergence results in  Theorem~\ref{th:drg} hold for \emph{any growth rate} of $c(N)$.
\end{remark}

\noindent
\textbf{Inhomogeneous levels of flexibility.}
Next, we consider a system that allows for inhomogeneous levels of flexibility for different task types. 
In particular, the compatibility graph is constructed by selecting each edge incident to a task-type $w \in W_N$ with probability $p_w(N)$, independently of other edges.
Thus, in expectation, task-type $w$ has the flexibility to be assigned to $Np_w(N)$ possible servers.
Theorem~\ref{th:erg} below provides a set of sufficient conditions on $(p_w(N))_{w\in W_N}$ to ensure that $\{G_N\}_{N \geq 1}$ satisfies the proportional sparsity and the subcriticality conditions.
\begin{theorem}[$\bigstar$]
\label{th:erg}
Assume that $(p_w(N))_{w\in W_N}$ satisfies the following:
\begin{equation}\label{eq:inhom-cond}
    \begin{aligned}
    {\normalfont (1)}\quad& \frac{N}{\ln(M(N)) + \ln(N)}\min_{w \in W_N} p_w(N) \to \infty,\\[10pt]
    {\normalfont (2)}\quad& M(N)\min_{w \in W_N} p_w(N)  \to \infty, \text{ and}\\[3pt]
    {\normalfont (3)}\quad& \frac{\ln(N)}{ (M(N))^2}\sum_{w\in W_N} \frac{1}{(p_w(N))^2} \to 0,
\end{aligned}
\end{equation}
as  $N\to\infty.$
Also, construct $G_N$ as follows: for any $v\in V_N, w\in W_N$, select edge $(v, w) \in E_N$ independently with probability $p_w(N)$.
Then the sequence $\{G_N\}_{N \geq 1}$ is proportionally sparse and satisfies the subcriticality condition, almost surely.
Consequently, the conclusions of Theorem~\ref{th:maintheorem} hold.
\end{theorem}
As in the proof of Theorem \ref{th:drg}, the proof of Theorem \ref{th:erg} also uses concentration of measure arguments, and it is provided in Appendix~\ref{app:erg}.

\begin{remark}\normalfont
If $p_w(N) = p(N)$ for all $w \in W_N$, then it is possible to relax condition (3) in equation \eqref{eq:inhom-cond} on the edge probabilities to $p(N) = \omega(\ln(N) / M(N))$.
For the special case $M(N) = N$, if $p(N) = o(\ln(N) / N)$ then the compatibility graph constructed as above, will leave at least one dispatcher isolated with probability tending one, as $N$ tends to infinity (see for example \cite{remco-book-1}). 
This means the graph sequence cannot satisfy the subcriticality condition, and this growth rate condition for $p(N)$ is nearly the optimum.
\end{remark}

\section{Proofs}
\label{sec:proofoverview}

To prove Theorem~\ref{th:maintheorem}, we will follow the usual interchange of limits argument.
In particular, the proof consists of three key steps. First, in Section~\ref{sec:overviewtransient}, we show that if the graph sequence is proportionally sparse, then the scaled occupancy process converges weakly to the appropriate system of ODEs, on any finite time interval.
Second, in Section~\ref{sec:stab-tight}, we show that if a graph sequence satisfies the subcriticality condition, then for any fixed $N$, the system is ergodic, and the sequence of steady state global occupancy is tight in the appropriate topology.
Third, in Section~\ref{sec:glob-stab}, we prove the global stability of the limiting system of ODEs with \eqref{eq:qinf} being its fixed point. 
Combining the above three steps we complete the proof of Theorem~\ref{th:maintheorem} in Section~\ref{sec:proof-main}.

\subsection{Process-level convergence}
\label{sec:overviewtransient}

In this section, we will prove the process-level convergence of the occupancy process under a general class of task assignment policies.
We start by specifying the class of assignment policies.\\

\noindent
\textbf{Task assignment policies.}
To determine which server to assign an incoming task to, each dispatcher in $W_N$ follows a generic task assignment policy $\Pi$ that works as follows. 
We identify the policy $\Pi$ with an \textit{assignment probability function} $\pp^\Pi = \left( p_0^\Pi, p_1^\Pi, \dots \right): \cX \to [0, 1]^\infty$. 
When a task arrives at a dispatcher $w \in W_N$ with LEQD $\xx^w$, 
\begin{enumerate}
    \item[(a)] Select a random queue length $I$ distributed as the probability measure induced by the assignment probability function evaluated at the LEQD $\xx^w$, that is, $\pro(I = i) =  p_i^\Pi(\xx^w)$ for $i=0,1,\ldots$.
    The random variable $I$ is independent of any other processes and also independent across different arrival epochs.
    \item[(b)] Next, a server is selected uniformly at random among all servers with queue length~$I$.
\end{enumerate}
Note that the above generic task assignment policy has a number of features, such as, the assignment policy should depend only on the LEQD, and the dispatcher does not distinguish between two neighboring servers with the same queue length.
It is not hard to see that the JSQ($d$) policy can also be described in the above form.
We restrict our analysis to Lipschitz continuous task assignment policies as given by the next definition.

\begin{definition}\textit{
\label{con:lip-routing}
(Lipschitz continuous task assignment policy)
A policy $\Pi$ is said to be Lipschitz continuous, if there exists a finite positive constant $K$, such that its assignment probability function satisfies the following. For any $\xx, \yy \in \cX$,
\begin{equation}\label{eq:lipschitz}
    \sum_{i=0}^\infty \left\lvert p_i^\Pi(\xx) - p_i^\Pi(\yy) \right\rvert \leq K  \sum_{i=0}^\infty \left\lvert x_i - y_i \right\rvert.
\end{equation}
}
\end{definition}

The Lipschitz continuity bounds the sensitivity of the assignment probability function.
In other words, if the $\ell_1$-distance between the LEQDs is small, then the probabilities of routing a task to a server with a given queue length from these states should also be close.
 As we will see in Lemma~\ref{lem:jsqdlipschitz}, the JSQ($d$) policy is Lipschitz continuous for any fixed $d\geq 1$.
 An example of a policy that is \emph{not} Lipschitz continuous is the ordinary JSQ policy, where the addition of a single task to the system, causing a change of $\Theta(1/N)$ in the LEQD, can change an assignment probability from 0 to 1.

\begin{lemma}[$\bigstar$]\label{lem:jsqdlipschitz}
For any fixed $d\geq 1$, the JSQ($d$) policy is Lipschitz continuous with Lipschitz constant $2d!\times d^2$.
\end{lemma}

The proof of Lemma \ref{lem:jsqdlipschitz} is provided in Appendix \ref{app:prooflemmajsqdlipschitz}. We now state the process-level convergence theorem.

\begin{theorem}
\label{th:transient}
Let $\{ G_N \}_{N \geq 1}$ be a proportionally sparse graph sequence and $\Pi$ be a Lipschitz continuous policy with 
assignment probability function $\pp^{\Pi} = \left( p_0^\Pi, p_1^\Pi, \dots \right)$.
Then, on any finite time interval $[0, T]$, the scaled occupancy process $\qq(\Phi(G_N, t)) = (q_1(\Phi(G_N, t)), q_2(\Phi(G_N, t)), \dots)$ converges weakly with respect to the Skorohod-$J_1$ topology, as $N\to\infty$, to the process $\qq^*(t) = (q_1^*(t), q_2^*(t), \dots)$, given by the unique solution of the system of ODEs:
\begin{equation}\label{eq:ODE-gentrans}
    \frac{d q_i^*(t)}{dt} = \lambda p_{i-1}^\Pi\big(( q_j^*(t) - q_{j+1}^*(t))_{ j \geq 0}\big) - (q_i^*(t) - q_{i+1}^*(t)) \quad\text{for}\quad i = 1, 2, \dots.
\end{equation}
provided $\qq^*(0) \in \cY$ and $\|\qq(\Phi(G_N, 0)) - \qq^*(0) \|_1 \to 0$ as $N \to \infty$.
\end{theorem}

The rest of this section is devoted to the proof of Theorem \ref{th:transient}.

\begin{remark}\normalfont
The unique solvability of the infinite set of ODE in equation \eqref{eq:ODE-gentrans} follows from the Lipschitz property of the policy $\Pi$, using standard results in analysis (see for example Theorem 3.2 in \cite{deimling06}).
\end{remark}

There are two main ingredients to the proof of Theorem \ref{th:transient}. 
First, in Proposition~\ref{lemma:queuedistr} we establish that for \emph{almost} all dispatchers, the LEQD is close to the GEQD, uniformly over any finite time interval. 
Second, we couple the $N$-th system with graph structure $G_N$ to a fully flexible system with complete bipartite graph and establish in Proposition~\ref{prop:qlessdelta}, a criterion for the two system to behave similarly. 
Finally, we complete the proof of Theorem \ref{th:transient}
using Propositions~\ref{lemma:queuedistr} and \ref{prop:qlessdelta}, and the Lipschitz continuity of the task assignment policy.

\subsubsection{Proximity of local and global empirical queue length distributions}
\label{sec:proximityqueuedistr}

We begin by introducing the notion of \textit{good} and \textit{bad} dispatchers. Loosely speaking, a \emph{good} dispatcher is a dispatcher for which the LEQD is close to the GEQD.

\begin{definition}
($\varepsilon$-good dispatchers)\label{def:e-good}
\textit{When the system is in state $z \in \cS_N$, the dispatcher $w \in W_N$ is called $\varepsilon$-good if
\begin{equation}
    \sum_{i=0}^\infty \lvert x_i(z) - x_i^w(z) \rvert < \varepsilon.
\end{equation}
Also, a dispatcher is called $\varepsilon$-bad if it is not $\varepsilon$-good.}
\end{definition}

In the following, let $B_N^\varepsilon(t)$ denote the number of $\varepsilon$-bad dispatchers at time $t$, when the system is in state $\Phi(G_N, t)$.

\begin{proposition}[$\bigstar$]
\label{lemma:queuedistr}
Let $\{ G_N \}_{N \geq 1}$ be a proportionally sparse graph sequence.  Then for each $\varepsilon, \delta > 0$,
\begin{equation}
\label{eq:badsetsmall}
    \mathbb{P}\Big( \sup_{t \in [0, T]} B_N^\varepsilon(t) \geq \delta M(N) \Big) \longrightarrow 0 \quad\text{as}\quad N \to \infty,
\end{equation}
provided $\qq^*(0) \in \cY$ and $\|\qq(\Phi(G_N, 0)) - \qq^*(0) \|_1 \to 0$ as $N \to \infty$.

\end{proposition}

The key idea in the proof of Proposition \ref{lemma:queuedistr} is to observe that the servers with queue length $i$ at time $t$ form a subset $U_i(t) \subseteq V_N$ and 
the proportional sparsity of the compatibility graph implies that \emph{for almost all} $w\in W_N$, the fraction of neighbors within $U_i(t)$ is close to $|U_i(t)|/N$.
However, we need to deal with some technical challenges, for example, the subset of dispatchers for which the above does \emph{not} hold, may depend on $i$, and thus, one needs to be careful in estimating the number of $\varepsilon$-bad dispatchers. Also, an uniformity over $[0,T]$ needs to be established.
The complete proof of Proposition \ref{lemma:queuedistr} is provided in Appendix \ref{app:queuedistr}. 

\begin{remark}\normalfont
It is worthwhile to highlight that Proposition~\ref{lemma:queuedistr} only requires the proportional sparsity property of the graph sequence and does \emph{not} depend on the task assignment policy or even on the dynamics in any way.
This makes the applicability of this method much broader, in analyzing certain structurally constrained large-scale dynamical systems where the process running at each vertex (queue length in our case) takes countably many values.
\end{remark}
 
\subsubsection{A coupling construction}
\label{sec:couplingconstruction}

Next, for any $N\geq 1$, we couple the queuing system on an arbitrary bipartite graph $G_N$ with the fully flexible system, that is, corresponding to the complete bipartite graph $K_{N,M}$. 

The coupling approach has been highly successful in proving large-system limit theorems. However, the biggest issue in constructing an appropriate coupling in the presence of compatibility constraints lies in the fact that if the arrivals are synchronized at each dispatcher $w \in W_N$, in two systems with different compatibility graphs, the set of neighbors of $w$ becomes different in the two systems. Consequently, one cannot synchronize which server a task will be assigned to, and the coupling of the queues breaks down. We will now introduce a novel coupling, called \textbf{optimal coupling}, to tackle this issue.

In short, we will refer to the two systems as $G_N$-system and $K_{N,M}$-system, respectively.
Both systems employ a Lipschitz continuous task assignment policy $\Pi$.
To describe the coupling, first, in each of the two systems, order the servers by non-decreasing queue lengths, breaking ties arbitrarily.\\

\noindent
$\bullet$ \textit{Departures}. Synchronize the departure epochs of the $k$-th ordered servers in the two systems, that is, both systems will potentially finish serving a task at the $k$-th ordered server at the same epoch, whenever they are non-empty for $k=1,2\ldots,N$.\\

\noindent
$\bullet$ \textit{Arrivals}. 
Synchronize the arrival epochs of task-type $w$ in both systems, for all $w\in W_N$.
At an arrival epoch of $w$, let $\xx^w$ and $\yy^w$ be the LEQDs for the two systems, respectively. 
Note that $\yy^w$ is also the GEQD for the $K_{N,M}$-system.
Define $p_i = \min(p_i^\Pi(\xx^w), p_i^\Pi(\yy^w))$ for $i=0,1,\ldots$. 
Now, let us draw a Uniform$[0, 1]$ random variable, independently of any other processes, and denote it by $U$.
Recall the description of the task assignment policy $\Pi$ from Section~\ref{sec:overviewtransient}.
The value of $U$ will be used in both systems to generate the random variables $I_1$ and $I_2$, for the $G_N$-system and $K_{N,M}$-system, respectively.
In the $G_N$-system, set $I_1 = i$, for $i=0,1,\ldots$, if
\begin{equation}
    U \in \Big[ \sum_{j = 0}^{i-1} p_j, \sum_{j = 0}^i p_j \Big) \cup \Big[ \sum_{j = 0}^\infty p_j + \sum_{j = 0}^{i-1} \big( p_j^\Pi(\xx^w) - p_j \big), \sum_{j = 0}^\infty p_j + \sum_{j = 0}^i \big( p_j^\Pi(\xx^w) - p_j \big) \Big),
\end{equation}
and assign the arriving task to a server selected uniformly at random among all servers with queue length $I_1$.
Similarly, in the $K_{N,M}$-system, set $I_2 = i$, for $i=0,1,\ldots$, if
\begin{equation}
    U \in \Big[ \sum_{j = 0}^{i-1} p_j, \sum_{j = 0}^i p_j \Big) \cup \Big[ \sum_{j = 0}^\infty p_j + \sum_{j = 0}^{i-1} \big( p_j^\Pi(\yy^w) - p_j \big), \sum_{j = 0}^\infty p_j + \sum_{j = 0}^i \big( p_j^\Pi(\yy^w) - p_j \big) \Big),
\end{equation}
and assign the arriving task to a server selected uniformly at random among all servers with queue length $I_2$.\\

Note that the coupling preserves the marginal laws of both systems. Next, we introduce a notion that facilitates 
comparison of the performance of the two systems on suitable asymptotic scales.

\begin{definition}\textit{
(Mismatch in queue length)
At an arrival epoch, the two coupled systems are said to mismatch in queue length, if $I_1 \neq I_2$, that is, the arriving tasks are assigned to two servers of different queue lengths in the two systems. 
Denote by $\Delta_N(t)$ the cumulative number of times the systems mismatch in queue length up to time $t$.}
\end{definition}
 The occupancy process and the cumulative number of mismatches in queue length are related as stated in Proposition~\ref{prop:qlessdelta} below.
The proof of Proposition~\ref{prop:qlessdelta} is provided in Appendix \ref{app:qlessdelta}.
\begin{proposition}[$\bigstar$]
\label{prop:qlessdelta}
For any $N\geq 1$, consider the $G_N$-system and the $K_{N,M}$-system coupled as above. Then the following holds almost surely on the coupled probability space: for all $t \geq 0$, 
\begin{equation}\label{eq:qlessdelta}
    \sum_{i=1}^\infty \left\lvert Q_i(\Phi(K_{N,M}, t)) - Q_i(\Phi(G_N, t)) \right\rvert \leq 2 \Delta_N(t),
\end{equation}
provided the inequality holds at $t=0$.
\end{proposition}

\subsubsection{Proof of Theorem \ref{th:transient}}
\label{sec:prooftransient}

\begin{proof}[Proof of Theorem \ref{th:transient}]
First, it is fairly standard (see for example, \cite[Chapter 8]{Kurtz81}) to show that $\qq(\Phi(K_{N,M}, t))$ converges weakly to $\qq^*(t)$.
Therefore, by Proposition~\ref{prop:qlessdelta}, it is enough to prove that for any $\varepsilon'>0$ and $\delta' > 0$ there exists $N_0 \geq 1$ such that
\begin{equation}\label{eq:local5.8}
\mathbb{P}(\sup_{t \in [0, T]} \Delta_N(t)/N \geq \varepsilon') < \delta' \text{ for all } N \geq N_0.
\end{equation}
Fix any $\varepsilon>0$, to be chosen later.
Recall Definition~\ref{def:e-good} and let $A_N^\varepsilon(t)$ and $B_N^\varepsilon(t)$ denote the number of $\varepsilon$-good and  $\varepsilon$-bad dispatchers in the $G_N$-system, respectively. Couple the $G_N$-system to the $K_{N,M}$-system by the optimal coupling described in Section~\ref{sec:couplingconstruction}. For brevity, let $\xx(t) := \xx(\Phi(G_N, t))$ and $\xx^w(t) := \xx^w(\Phi(G_N, t))$ for the $G_N$-system and $\yy(t) := \xx(\Phi(K_{N,M}, t))$ for the $K_{N,M}$-system. Also, define $\rho_N(t) := \sum_{i=0}^\infty \lvert y_i(t) - x_i(t) \rvert$.
At an arrival epoch $t \geq 0$, if a task arrives at an $\varepsilon$-good dispatcher $w \in W_N$ then
\begin{equation}
\begin{gathered}
    \sum_{i=0}^\infty \lvert y_i(t-) - x_i^w(t-) \rvert \leq \sum_{i=0}^\infty \lvert y_i(t-) - x_i(t-) \rvert + \sum_{i=0}^\infty \lvert x_i(t-) - x_i^w(t-) \rvert 
    \leq \rho_N(t-) + \varepsilon.
\end{gathered}
\end{equation}
Define the uniform random variable $U$ and $p_i$ as in the optimal coupling, and observe that the probability that the systems mismatch in queue length at such an arrival epoch, is bounded by
\begin{equation}
\begin{split}
    &\mathbb{P}\Big( U \notin \big[ 0, \sum_{i = 0}^\infty p_i \big) \Big) = 1 - \sum_{i = 0}^\infty p_i = \sum_{i = 0}^\infty \big( p_i^\Pi(\yy(t-)) - p_i \big) \\
    &\leq \sum_{i = 0}^\infty \left\lvert p_i^\Pi(\yy(t-)) - p_i^\Pi(\xx^w(t-)) \right\rvert
    \leq K \sum_{i=0}^\infty \left\lvert y_i(t-) -  x_i^w(t-) \right\rvert \leq K \left( \rho_N(t-) + \varepsilon \right),
\end{split}
\end{equation}
by the Lipschitz property of $\Pi$ in~\eqref{eq:lipschitz}. 
If instead, a task arrives at an $\varepsilon$-bad dispatcher, then with probability at most one the systems mismatch in queue length. 
Now, since at the arrival epochs, the random variables $U$ are independent of any other processes, we can construct an independent unit-rate Poisson process $(Z(t))_{t\geq 0}$, so that $\Delta_N(t)$ is upper bounded by a random time-change (cf.~\cite[\S 2.1]{PTRW07}) of $Z$ as follows:
for all $t\in [0,T]$,
\begin{equation}\label{eq:local5.11}
\begin{gathered}
    \Delta_N(t) \leq Z\Big( \frac{\lambda N}{M(N)} \int_0^t \big[A_N^\varepsilon(s-) \cdot K (\rho_N(s-) + \varepsilon) + B_N^\varepsilon(s-) \cdot 1 \big]\diff s \Big).
\end{gathered}
\end{equation}
Now observe that
\begin{equation}\label{eq:local5.13}
\begin{split}
    \rho_N(t) &= \sum_{i=0}^\infty \left\lvert x_i(\Phi(K_{N,M}, t) - x_i(\Phi(G_N, t)) \right\rvert \\
   & = \sum_{i=0}^\infty \left\lvert q_i(\Phi(K_{N,M}, t)) - q_{i+1}(\Phi(K_{N,M}, t)) - q_i(\Phi(G_N, t)) + q_{i+1}(\Phi(G_N, t)) \right\rvert \\
    &\leq 2 \sum_{i=0}^\infty \left\lvert q_i(\Phi(K_{N,M}, t)) - q_i(\Phi(G_N, t)) \right\rvert \leq \frac{4 \Delta_N(t)}{N},
\end{split}
\end{equation}
where the last inequality follows from Proposition \ref{prop:qlessdelta}.
Furthermore, using Proposition \ref{lemma:queuedistr} to bound $B_N^\varepsilon$ on the right hand side of~\eqref{eq:local5.11}, we can write, for any fixed $\delta>0$ to be chosen later,
there exists $N_0=N_0(\delta)\in \N$, such that for all $N\geq N_0$,
\begin{equation}\label{eq:local5.12}
    \mathbb{P}\Big(\sup_{t\in [0, T]} B_N^\varepsilon(t)> \delta M(N) \Big) \leq \frac{\delta}{2}.
\end{equation}
Therefore,~\eqref{eq:local5.13},~\eqref{eq:local5.12}, and an application of Tonelli's theorem imply that for all $N\geq N_0$ and $t\in [0,T]$,
\begin{equation}\label{eq:local5.14}
\begin{split}
    \mathbb{E}\Big( \frac{\Delta_N(t)}{N} \Big)
    &\leq \lambda \int_0^t \Big[K \Big( 4 \mathbb{E}\Big( \frac{\Delta_N(s-)}{N} \Big) + \varepsilon \Big) + \frac{3\delta}{2} \Big] \diff s.
\end{split}
\end{equation}
By applying Gr\"{o}nwall's inequality to~\eqref{eq:local5.14}, we obtain,
\begin{equation}
    \mathbb{E}\Big( \frac{\Delta_N(t)}{N} \Big) \leq \lambda (K \varepsilon + 3\delta/2) t \exp(4 \lambda K t).
\end{equation}
Finally, since $\Delta_N(t)$ is nonnegative, and nondecreasing in $t$, using Markov's inequality, we write,
\begin{equation}
    \begin{split}
    \mathbb{P}\Big(\sup_{t\in [0,T]} \frac{\Delta_N(t)}{N} \geq \varepsilon' \Big)
    &\leq \mathbb{P}\Big( \frac{\Delta_N(T)}{N} \geq \varepsilon' \Big)
    \leq \frac{1}{\varepsilon'} \mathbb{E}\Big( \frac{\Delta_N(T)}{N} \Big) \\
    &\leq \frac{\lambda}{\varepsilon'} (K \varepsilon + 3\delta/2) T \exp(4 \lambda K T),
    \end{split}
\end{equation}
which yields equation \eqref{eq:local5.8} by choosing $\varepsilon$ and $\delta$ small enough.
\end{proof}

\subsection{Stability and tightness}\label{sec:stab-tight}

In this section we will prove positive recurrence of the Markov process $\Phi(G_N, \cdot)$ and tightness of the steady state occupancy process, as stated in the next theorem.
\begin{theorem}
\label{th:stability}
Let $\{ G_N \}_{ N \geq 1 }$ be a sequence of graphs satisfying the subcriticality condition. Then
\begin{enumerate}[label=(\roman*)]
\item There exists $N_0 \geq 1$, such that for each $N \geq N_0$, the Markov process $\Phi(G_N, \cdot)$ under the JSQ($d$) policy with $d\geq 1$ is positive Harris recurrent, and thus, it has a unique steady state.
\item Let $\Phi(G_N, \infty)$ denote a random variable distributed as the steady state of $\Phi(G_N, \cdot)$. 
Fix any $r \in (1, 2 / (1 + \lambda))$ and a positive sequence $\boldsymbol{\omega} = (\omega_1, \omega_2,\ldots)$ satisfying
\begin{equation}\label{eq:omega-cond}
    \omega_{i_0 + i} = \omega_{i_0} r^i \text{ for } i \geq 1
\end{equation}
for some $i_0 \in \mathbb{N}$. 
Then the sequence $\{ \qq(\Phi(G_N, \infty)) \}_{N \geq 1}$ is tight when $\cY$ is endowed with the $\ell_1^\omega$-topology.
\end{enumerate}
\end{theorem}

The positive Harris recurrence from part (i) of Theorem~\ref{th:stability} follows from known techniques in \cite{Bramson11} or \cite{CBL19}. However, as we prove the tightness in part (ii), the positive recurrence follows as a side result. This is given in Appendix~\ref{app:proofstability}. In the rest of this section, we complete the proof of tightness. First, we obtain a bound on the tail of the expected global occupancy of the stationary state.

\label{sec:tightness}


\begin{lemma}[$\bigstar$]
\label{lemma:momentbound}
Consider a sequence $\{ G_N \}_{N \geq 1}$ that satisfies the subcriticality condition. There exists $N_0 \geq 1$ such that for all $N \geq N_0$ and $k \geq 1$,
\begin{equation}
\label{eq:constraintocc}
     \sum_{i = k}^\infty \mathbb{E}[q_i(\Phi(G_N, \infty))]  \leq \frac{(1 + \lambda) / 2}{1 - (1 + \lambda) / 2} \mathbb{E}[q_{k - 1}(\Phi(G_N, \infty))].
\end{equation}
\end{lemma}
The proof of Lemma \ref{lemma:momentbound} is provided in Appendix \ref{app:proofstability}.
It uses a sequence $(V_k)_{k\geq 1}$ of Lyapunov functions, where bounds on the drift of $V_k$ result in a moment bound on the tail-sum of $\qq(\Phi(G_N, \infty))$ starting from $k$. 
The next technical lemma states sufficient criteria for tightness with respect to the $\ell_1^\omega$-topology, and is proved in Appendix~\ref{app:tightnessequiv}. 

\begin{lemma}[$\bigstar$]
\label{lemma:tightnessequiv}
Fix a positive sequence $\boldsymbol{\omega} = (\omega_1, \omega_2,\ldots)$.
The sequence of random variables $\{ \qq(\Phi(G_N, \infty)) \}_{N \geq 1}$ is tight when $\cY$ is endowed with the $\ell_1^\omega$-topology, if  for each fixed $i\in \N$, $\{ q_i(\Phi(G_N, \infty)) \}_{N \geq 1}$ is tight in $\R$ and
\begin{equation}
\label{eq:tightvanishingtail}
    \lim_{j \rightarrow \infty} \limsup_{N \rightarrow \infty} \mathbb{P}\left( \sum_{i = j}^\infty \omega_i q_i(\Phi(G_N, \infty)) > \varepsilon \right) = 0.
\end{equation}
\end{lemma}

\begin{proof}[Proof of Theorem \ref{th:stability}(ii)]
We will verify the sufficient conditions stated in Lemma~\ref{lemma:tightnessequiv}.
Because $\omega_i q_i(\Phi(G_N, \infty)) \in [0, \omega_i]$, it is trivially tight in $\R$. Now, by Markov's inequality,
\begin{equation}
\label{eq:markovocc}
    \mathbb{P} \Big( \sum_{i = j}^\infty \omega_i q_i(\Phi(G_N, \infty)) > \varepsilon \Big) \leq \frac{1}{\varepsilon} \mathbb{E}\Big[ \sum_{i = j}^\infty \omega_i q_i(\Phi(G_N, \infty)) \Big] = \frac{1}{\varepsilon} \sum_{i = j}^\infty \omega_i \mathbb{E}\left[ q_i(\Phi(G_N, \infty)) \right].
\end{equation}
Fix any $r \in (1, 2 / (1 + \lambda))$.
In the following, let $j \geq i_0$ be such that $\omega_i = (\omega_{i_0} / r^{i_0}) r^i$ for $i \geq j$. 
Recall from Lemma~\ref{lemma:momentbound} that $\{ \mathbb{E}[q_i(\Phi(G_N, \infty))] \}_{i\geq 1}$ satisfies~\eqref{eq:constraintocc}.
We will bound the right-hand side of equation \eqref{eq:markovocc}, by maximizing its value over all sequences in $\ell_1$ satisfying equation~\eqref{eq:constraintocc}. That is, we arrive at the following maximization problem: for $j=1,2,\ldots$,
\begin{equation}
\begin{gathered}
    \textsc{Sol}(j) := (\omega_{i_0} / r^{i_0}) \sup_{\boldsymbol{a} \in \ell_1} \sum_{i = j}^\infty r^i a_i, \text{ where } 
    \sum_{i = k}^\infty a_i \leq \frac{(1 + \lambda) / 2}{1 - (1 + \lambda) / 2} a_{k - 1} \text{ for all } k \geq 1,
\end{gathered}
\end{equation}
and its finite projection:
\begin{equation}
\begin{gathered}
    \textsc{Sol}_n(j) = (\omega_{i_0} / r^{i_0}) \sup_{\boldsymbol{a} \in \mathbb{R}^n} \sum_{i = j}^n r^i a_i, \text{ where } 
    \sum_{i = k}^n a_i \leq \frac{(1 + \lambda) / 2}{1 - (1 + \lambda) / 2} a_{k - 1} \text{ for all } 1 \leq k \leq n,
\end{gathered}
\end{equation}
Note that the solution of the maximization problem is an upper bound on the right-hand side of \eqref{eq:markovocc}, and thus
\begin{equation}
    \sum_{i = j}^\infty \omega_i \mathbb{E}\left[ q_i(\Phi(G_N, \infty)) \right] \leq \textsc{Sol}(j) = \sup_{n \geq 1} \textsc{Sol}_n(j),
\end{equation}
where the last equality follows because each sequence in $\ell_1$ can be written as the limit of a sequence of vectors in $\mathbb{R}^n$ in the $\ell_1$-topology. The finite projected maximization problem is linear and hence the solution is found at one of the boundary points of the feasible region. In this case, for fixed $n \geq 1$ and $j \geq i_0$, there is only one boundary point, the point for which the constraints hold with equality:
\begin{equation}
\begin{aligned}
    a_i = \left( \frac{1 + \lambda}{2} \right)^i \text{ for } 0 \leq i \leq n - 1, &&
    a_n = \frac{1}{1 - (1 + \lambda) / 2} \left( \frac{1 + \lambda}{2} \right)^n,
\end{aligned}
\end{equation}
such that
\begin{equation}
    \textsc{Sol}_n(j) = (\omega_{i_0} / r^{i_0}) \sup_{a \in \mathbb{R}^n} \sum_{i = j}^n r^i a_i = \frac{\omega_{i_0}}{r^{i_0}} \frac{r^j}{1 - (1 + \lambda) / 2} \left( \frac{1 + \lambda}{2} \right)^j.
\end{equation}
Therefore,
\begin{equation}
    \lim_{j \rightarrow \infty} \limsup_{N \rightarrow \infty} \mathbb{P}\left( \sum_{i = j}^\infty \omega_i q_i(\Phi(G_N, \infty)) > \varepsilon \right) \leq \lim_{j \rightarrow \infty} \frac{1}{\varepsilon} \frac{\omega_{i_0}}{r^{i_0}} \frac{r^j}{1 - (1 + \lambda) / 2} \left( \frac{1 + \lambda}{2} \right)^j = 0.
\end{equation}

\end{proof}

\subsection{Global stability}\label{sec:glob-stab}

The last key property needed in establishing the interchange of limits is the global stability of the process-level limit in~\eqref{eq:ODE-gen}. As stated in the next theorem, global stability shows that the limiting system of ODEs converges to a fixed point if started from suitable states.

\begin{theorem}
\label{th:globalstability}
For any positive sequence $\boldsymbol{\omega} = (\omega_1, \omega_2,\ldots)$, let us denote 
$$\Psi_\omega(t) := \sum_{i = 1}^\infty \omega_i \lvert q^*_i(t) - q_i^*(\infty) \rvert$$ 
and $i_0 := \min\{ i \geq 1 \mid \lambda (2 q_i^*(\infty) + 1) < \frac{1 + \lambda}{2} \}$. For each $\lambda < 1$, there exists a choice of $\boldsymbol{\omega} = (\omega_1, \omega_2,\ldots)$, satisfying
\begin{equation}
\label{eq:globalstabomega}
    \omega_{i_0 + i} = \omega_{i_0} r^i \text{ for } i \geq 1
\end{equation}
for some $r \in (1, 2 / (1 + \lambda))$, such that $\Psi_\omega(t)$ converges exponentially to zero as $t\to\infty$, if $\Psi_\omega(0) < \infty$.
\end{theorem}
The proof relies on the global stability result of the classical meanfield limit of the JSQ($d$) policy, as stated in \cite[Theorem 3.6]{Mitzenmacher01}. The details can be found in Appendix \ref{app:proofglobalstability}.

\subsection{Proof of Theorem \ref{th:maintheorem}}\label{sec:proof-main}

We now have all the necessary results to prove Theorem \ref{th:maintheorem}.

\begin{proof}[Proof of Theorem \ref{th:maintheorem}]
The process-level convergence of the occupancy process under proportional sparsity, follows from Theorem \ref{th:transient}. 
For the convergence of steady states, we follow the usual interchange of limits argument.

Fix $\boldsymbol{\omega} = (\omega_1, \omega_2,\ldots)$ as in equation \eqref{eq:globalstabomega} for $r$ chosen as in Theorem \ref{th:globalstability}.
By Theorem \ref{th:stability}, $\{ \qq(\Phi(G_N, \infty)) \}_{N \geq 1}$ is tight when $\cY$ is endowed with the $\ell_1^\omega$-topology and hence, any subsequence has a convergent further subsequence. Let $\{ \qq(\Phi(G_{N_n}, \infty)) \}_{n \geq 1}$ be such a convergent subsequence with $\{ N_n \}_{n \geq 1} \subseteq \mathbb{N}$ and assume that $\qq(\Phi(G_{N_n}, \infty))$ converges weakly to some distribution $\hat{\pi}$, as $n \rightarrow \infty$. 
Observe that $\hat{\pi}$ must be supported on
$$\cY_{\omega} := \Big\{ \qq \in \ell_1 \ \Big| \sum_{i=1}^\infty \omega_i q_i <\infty, \ q_i\in [0,1], q_i\geq q_j, i<j, i,j\in \mathbb{N} \Big\}.$$
Now, imagine starting the system in the steady state as $\Phi(G_{N_n}, 0) \sim \Phi(G_{N_n}, \infty)$. Then $\Phi(G_{N_n}, t) \sim \Phi(G_{N_n}, \infty)$ for all $t \in [0, T]$, and also, 
by Theorem \ref{th:transient}, it follows that
\begin{equation}
    \qq(\Phi(G_{N_n}, t)) \xrightarrow{d} \qq^*(t) \text{ as } n \rightarrow \infty,
\end{equation}
uniformly on $[0,T]$.
 Therefore,
\begin{equation}
    \qq^*(t) \overset{d}{=} \hat{\pi} \text{ for all } t \in [0, T].
\end{equation}
It follows that $\hat{\pi}$ is an invariant distribution of the limiting dynamics. To this end, the global stability result in Theorem~\ref{th:globalstability} implies that $\hat{\pi}$ must be the Dirac measure at the fixed point $\qq^*(\infty)$ of $\qq^*(t)$. 
This completes the proof of interchange of limits.
\end{proof}
\section{Simulation experiments}
\label{sec:simulation}

In this section, we perform extensive simulation experiments to evaluate our results and to gain insights into the cases that are not included in our theoretical framework, for example, when the compatibility graph may not be proportionally sparse. \\

\begin{figure}
\centerline{\begin{subfigure}[t]{.5\textwidth}
    \centering
    \includegraphics[width=\textwidth]{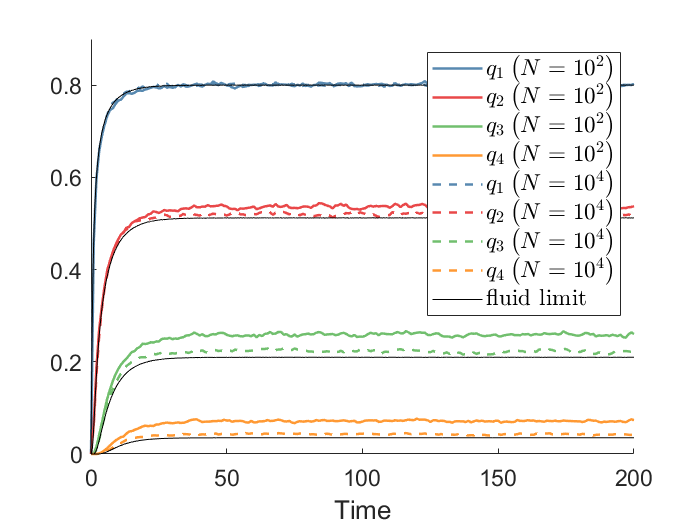}
    \caption{The occupancy process for compatibility graphs with average degree $(\ln(N))^2$ for $N = 10^2, 10^4$.}
    \label{fig:erg_clique_log2}
\end{subfigure}\hspace{2mm}%
\begin{subfigure}[t]{.5\textwidth}
    \centering
    \includegraphics[width=\textwidth]{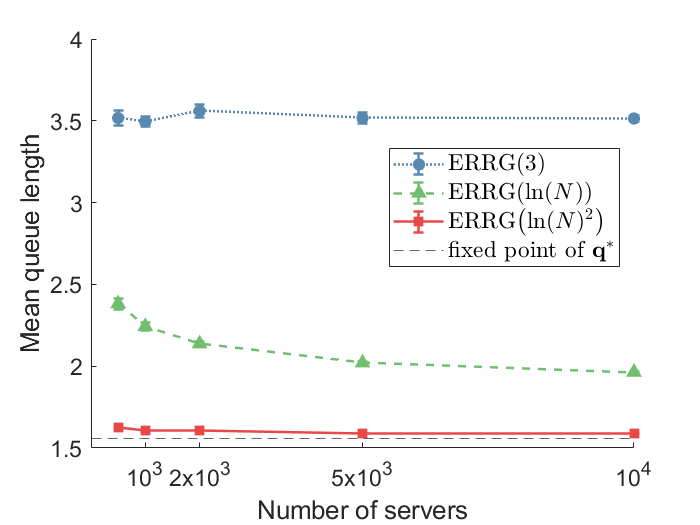}
    \caption{The average queue length in steady state for varying degrees.}
    \label{fig:erg_types}
\end{subfigure}
}
\caption{Performance of the JSQ($2$) policy on a bipartite Erd\H{o}s-R\'{e}nyi graph for $\lambda = 0.8$.}
\end{figure}

\noindent
\textbf{Verification of process-level convergence.}
Figure \ref{fig:erg_clique_log2} shows the evolution of the occupancy process for $M(N) = N = 10^2$ and $M(N) = N = 10^4$, where the compatibility graphs are single instances of bipartite Erd\H{o}s-R\'{e}nyi random graphs (ERRG) with edge probability $(\ln(N))^2/N$. Although the average number of neighbors is significantly less than the number of neighbors in a fully flexible system, the simulation illustrates that the occupancy processes closely follow the limiting dynamics of a fully flexible system. The sample path trajectories for $N=10^2$ further exhibits the validity of the asymptotic results for fairly small values of $N$.\\

\noindent
\textbf{Influence of the level of connectivity.}
In Figure \ref{fig:erg_types}, we examine the dependence of the average queue length in steady state on the average server-degree in $G_N$. In each case, we have assumed $M(N)=N$ and the compatibility graph is a single instance of the bipartite ERRG.
When the average degree scales as $(\ln(N))^2$, the graph sequence satisfies the conditions in Theorem~\ref{th:maintheorem}, and the average queue length can be seen to converge to the fixed point of $\qq^*(\infty)$ given by~\eqref{eq:qinf}. 
Among the graph sequences that do not satisfy the conditions in Theorem~\ref{th:maintheorem}, in the cases when the average degree is (i) a constant and (ii) equal to $\ln(N)$, the simulation shows that the behavior differs significantly from the fully flexible case.
 \\

\begin{figure}
\centerline{
\begin{subfigure}[t]{.5\textwidth}
    \centering
    \includegraphics[width=\textwidth]{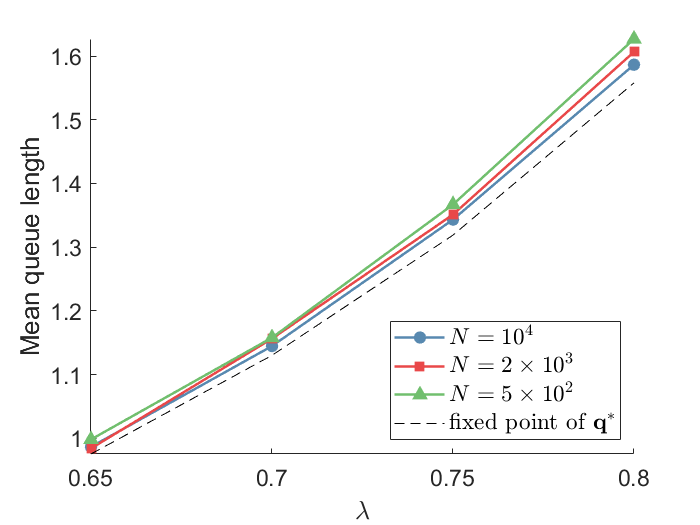}
    \caption{The average queue length in steady state for $(\ln(N))^2$ average degree and varying arrival rates.}
    \label{fig:erg_lambda}
\end{subfigure}\hspace{2mm}%
\begin{subfigure}[t]{.5\textwidth}
    \centering
    \includegraphics[width=\textwidth]{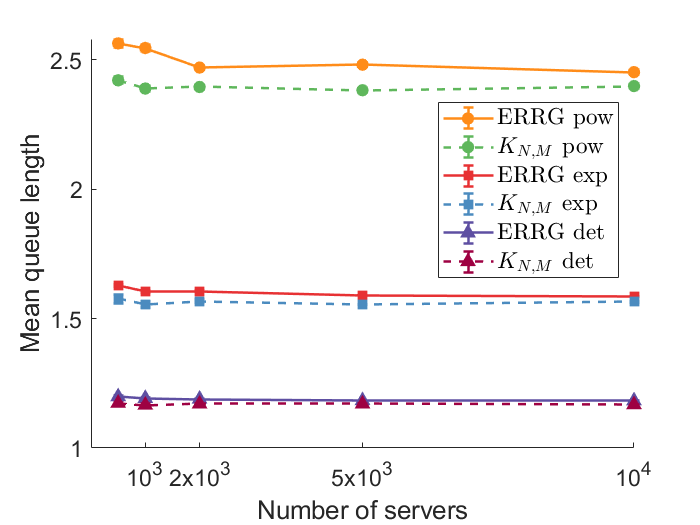}
    \caption{The average queue length in steady state for $(\ln(N))^2$ average degree and varying service distributions.}
    \label{fig:erg_det_W}
\end{subfigure}
}
\caption{Performance of the JSQ($2$) policy on a bipartite Erd\H{o}s-R\'{e}nyi random graph.}
\end{figure}

\noindent
\textbf{Influence of the arrival rate.}
Figure \ref{fig:erg_lambda} shows the dependence of the average queue length in steady state on the arrival rate $\lambda$. Although the average queue length converges to the fixed point $\qq^*(\infty)$ for all $\lambda < 1$, the rate of convergence is faster for smaller arrival rates. The difference in convergence rate is especially notable for $\lambda = 0.8$. \\

\noindent
\textbf{Influence of the service-time distribution.}
Figure \ref{fig:erg_det_W} shows the average queue length in steady state when service times are either exponentially distributed, deterministic -- fixed at one, or distributed as a power law with exponent 3. Note that these simulations go beyond the theoretical results for the exponential distribution. 
In each of the above three cases, Figure~\ref{fig:erg_det_W} shows that the behavior of the steady state average queue length with sparse bipartite ERRG compatibility graph is close to the corresponding fully flexible system.
This hints at a possible universality result, where one would expect to see an analog of Theorem~\ref{th:maintheorem} under a general service time distribution, although the behavior for different service time distributions will be different.\\

\begin{figure}
\centerline{\begin{subfigure}[t]{.35\textwidth}
    \centering
    {%
    \setlength{\fboxsep}{1pt}%
    \setlength{\fboxrule}{0.5pt}%
    \fbox{\includegraphics[width=\textwidth]{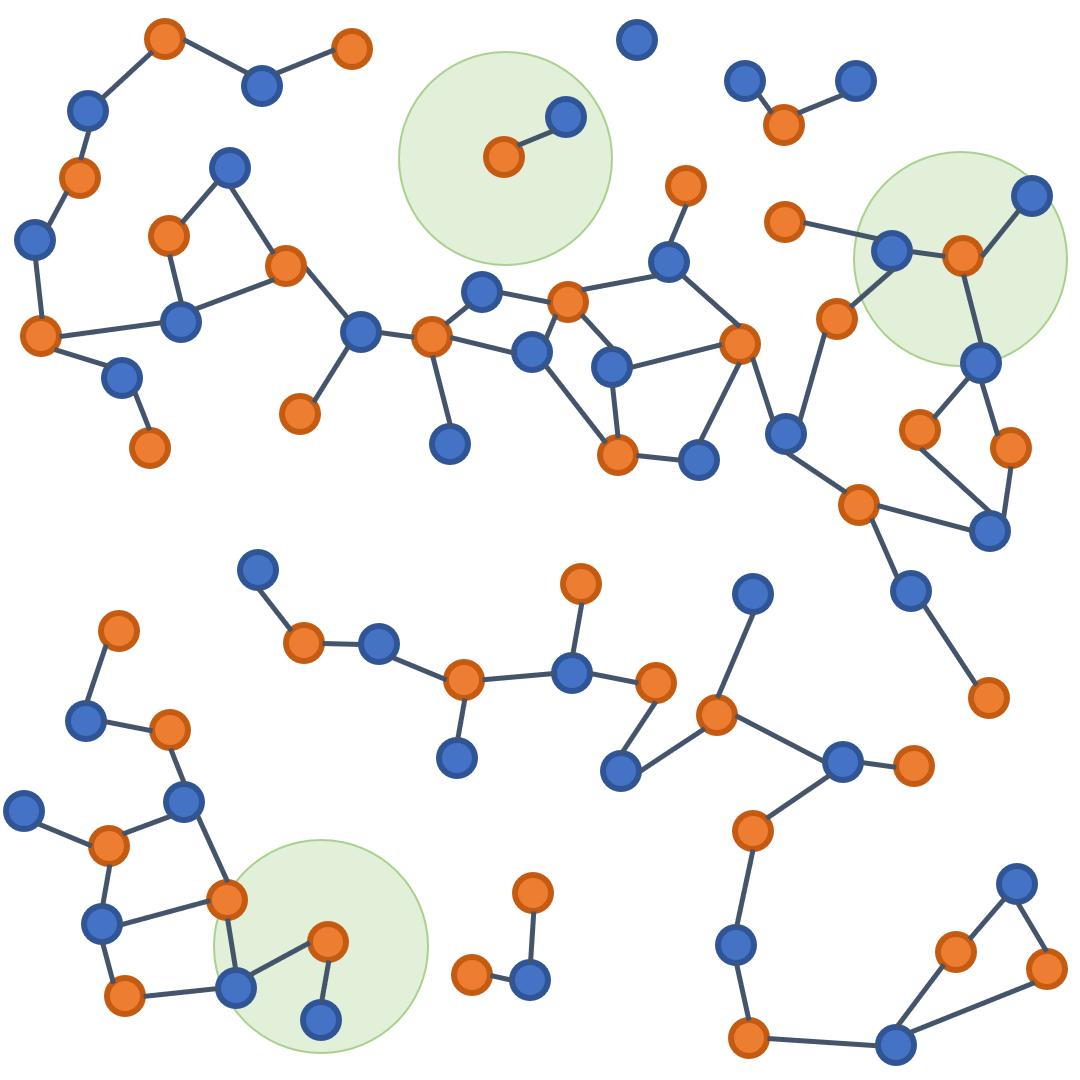}}%
    }%
    \caption{A random geometric graph is built by connecting dispatchers (in orange) to servers (in blue) that are within a specific distance to each other.}
    \label{fig:rggcircle}
\end{subfigure}\hspace{5mm}%
\begin{subfigure}[t]{.5\textwidth}
    \centering
    \includegraphics[width=\textwidth]{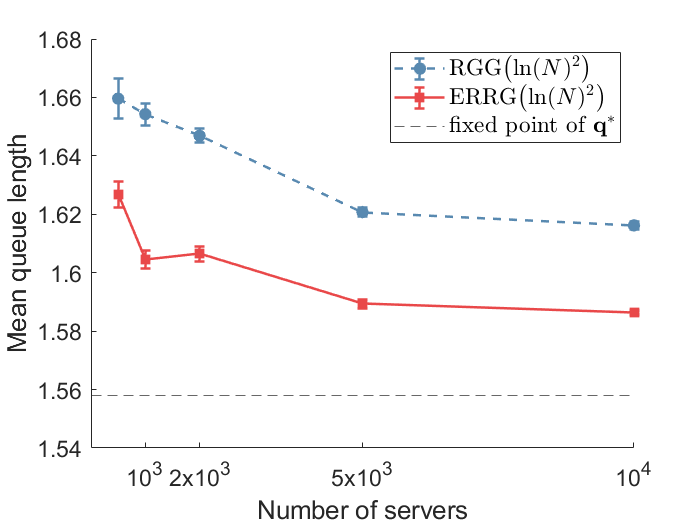}
    \caption{The average queue length in steady state of a bipartite Erd\H{o}s-R\'{e}nyi random graph for $(\ln(N))^2$ average degree, $d = 2$, and $\lambda = 0.8$, in comparison to a random geometric graph with same average degree.}
    \label{fig:erg_rgg}
\end{subfigure}}
\caption{Performance of large-scale systems with a spatially constrained compatibility graph.}
\end{figure}

\noindent
\textbf{Influence of geometry.}
In the broader context of spatial queueing systems, servers and dispatchers are often geometrically constrained. 
The compatibility graph in this case arises due to the proximity of servers and dispatchers and such a geometrically constrained graph often does not satisfy the proportional sparsity condition. 
In Figure~\ref{fig:rggcircle}, we investigate using the well-known \emph{random geometric graph} model, if and when the performance of such networks is asymptotically indistinguishable from the fully flexible model. 
A random geometric graph is built by placing each of the $N$ servers and $M(N)$ dispatchers at a uniformly selected location in $[0, 1]^2$.
As Figure \ref{fig:rggcircle} illustrates, dispatchers are then connected to all servers within a certain radius. Dispatchers will therefore have only local connections to servers in their proximity.

Figure \ref{fig:erg_rgg} shows a comparison between the average steady state queue length for two systems with compatibility graphs that are single instances of a bipartite ERRG and a random geometric graph, respectively. Both graphs have the same average server degree of $(\ln(N))^2$. 
Notice that the average queue length of the random geometric graph converges to the fixed point $\qq^*(\infty)$, although the proportional sparsity condition is not satisfied. Extending our results for such spatial queueing systems would be an important future research direction.

\section{Conclusion}
\label{sec:conclusion}

In this paper, we studied the impact of task-server compatibility constraints on the performance of the JSQ($d$) algorithm in large-scale systems.
The results extend the validity of the meanfield approximation far beyond the fully flexible scenario, allowing room for a much sparser class of compatibility graphs to asymptotically match the performance of a fully flexible system.
We have also provided explicit constructions of random compatibility graphs that are shown to almost surely belong to the above class, while being significantly sparser than the complete bipartite graph.
In the context of large-scale data centers, this translates into significant reductions in implementation complexity and storage capacity, since both are roughly proportional to the server-degrees in the compatibility graph.
Extensive simulation experiments corroborate the theoretical results.
More interestingly, some simulation experiments seem to suggest that the meanfield approximation may be valid for an even larger class of spatial graphs, where our expansion criterion does not hold.
This would be an interesting future research direction.
Also, in an ongoing work, we are exploring the heterogeneous case where the arrival and service rates may depend on the task type and the server, respectively,  to see how the sufficient condition on the graph sequence depends on the heterogeneity in the arrival rates or the service rates.

\bibliographystyle{apa}
\bibliography{main,mendeley_debankur}

\begin{thebibliography}{}

\bibitem[\protect\astroncite{Arapostathis et~al.}{2012}]{AHP12}
Arapostathis, A., Hmedi, H., and Pang, G. (2012).
\newblock {On uniform exponential ergodicity of Markovian multiclass
  many-server queues in the Halfin-Whitt regime}.
\newblock {\em Queueing Syst.}, 71:25--51.

\bibitem[\protect\astroncite{Braess}{1968}]{Braess68}
Braess, D. (1968).
\newblock {{\"{U}}ber ein Paradoxon aus der Verkehrsplanung}.
\newblock {\em Unternehmensforschung}, 12(1):258--268.

\bibitem[\protect\astroncite{Bramson}{2011}]{Bramson11}
Bramson, M. (2011).
\newblock {Stability of join the shortest queue networks}.
\newblock {\em Ann. Appl. Probab.}, 21(4):1568--1625.

\bibitem[\protect\astroncite{Budhiraja et~al.}{2019}]{BMW17}
Budhiraja, A., Mukherjee, D., and Wu, R. (2019).
\newblock {Supermarket model on graphs}.
\newblock {\em Ann. Appl. Probab.}, 29(3):1740--1777.

\bibitem[\protect\astroncite{Cardinaels et~al.}{2020}]{CBL20}
Cardinaels, E., Borst, S., and van Leeuwaarden, J. S.~H. (2020).
\newblock Redundancy scheduling with locally stable compatibility graphs.

\bibitem[\protect\astroncite{Cardinaels et~al.}{2019}]{CBL19}
Cardinaels, E., Borst, S.~C., and van Leeuwaarden, J. S.~H. (2019).
\newblock {Job assignment in large-scale service systems with affinity
  relations}.
\newblock {\em Queueing Syst.}, 93(3-4):227--268.

\bibitem[\protect\astroncite{Chung et~al.}{1989}]{CGW89}
Chung, F. R.~K., Graham, R.~L., and Wilson, R.~M. (1989).
\newblock Quasi-random graphs.
\newblock {\em Combinatorica}, 9(4):345--362.

\bibitem[\protect\astroncite{Cruise et~al.}{2020}]{CJS20}
Cruise, J., Jonckheere, M., and Shneer, S. (2020).
\newblock Stability of jsq in queues with general server-job class
  compatibilities.

\bibitem[\protect\astroncite{Deimling}{2006}]{deimling06}
Deimling, K. (2006).
\newblock {\em Ordinary differential equations in Banach spaces}, volume 596.
\newblock Springer.

\bibitem[\protect\astroncite{Ethier and Kurtz}{2009}]{EK2009}
Ethier, S.~N. and Kurtz, T.~G. (2009).
\newblock {\em {Markov Processes: Characterization and Convergence}}.
\newblock John Wiley {\&} Sons.

\bibitem[\protect\astroncite{Foss and Chernova}{1998}]{FC98}
Foss, S.~G. and Chernova, N.~I. (1998).
\newblock On the stability of a partially accessible multi-station queue with
  state-dependent routing.
\newblock {\em Queueing Syst.}, 29(1):55--73.

\bibitem[\protect\astroncite{Gamarnik and Stolyar}{2012}]{GS12}
Gamarnik, D. and Stolyar, A.~L. (2012).
\newblock Multiclass multiserver queueing system in the halfin-whitt heavy
  traffic regime: asymptotics of the stationary distribution.
\newblock {\em Queueing Syst.}, 71(1-2):25--51.

\bibitem[\protect\astroncite{Gast}{2015}]{Gast15}
Gast, N. (2015).
\newblock {The power of two choices on graphs: the pair-approximation is
  accurate}.
\newblock In {\em Proc. MAMA workshop 2015}, pages 69--71.

\bibitem[\protect\astroncite{Hajek}{1982}]{Hajek82}
Hajek, B. (1982).
\newblock {Hitting-time and occupation-time bounds implied by drift analysis
  with applications}.
\newblock {\em Adv. Appl. Prob.}, 14(3):502--525.

\bibitem[\protect\astroncite{Hmedi et~al.}{2019}]{HAP19}
Hmedi, H., Arapostathis, A., and Pang, G. (2019).
\newblock Uniform stability of a class of large-scale parallel server networks.

\bibitem[\protect\astroncite{Janson et~al.}{2000}]{JLR00}
Janson, S., Luczak, T., and Rucinski, A. (2000).
\newblock {\em {Random graphs}}.
\newblock John Wiley {\&} Sons.

\bibitem[\protect\astroncite{Kurtz}{1981}]{Kurtz81}
Kurtz, T.~G. (1981).
\newblock {\em Approximation of population processes}.
\newblock SIAM.

\bibitem[\protect\astroncite{Meyn and Tweedie}{1993}]{MT93}
Meyn, S.~P. and Tweedie, R.~L. (1993).
\newblock {\em {Markov Chains and Stochastic Stability}}.
\newblock Springer London.

\bibitem[\protect\astroncite{Mishra et~al.}{2010}]{MHCD10}
Mishra, A.~K., Hellerstein, J.~L., Cirne, W., and Das, C.~R. (2010).
\newblock Towards characterizing cloud backend workloads: insights from google
  compute clusters.
\newblock {\em ACM SIGMETRICS Performance Evaluation Review}, 37(4):34--41.

\bibitem[\protect\astroncite{Mitzenmacher}{1996}]{Mitzenmacher96}
Mitzenmacher, M. (1996).
\newblock {\em {The power of two choices in randomized load balancing}}.
\newblock PhD thesis, University of California, Berkeley.

\bibitem[\protect\astroncite{Mitzenmacher}{2001}]{Mitzenmacher01}
Mitzenmacher, M. (2001).
\newblock {The power of two choices in randomized load balancing}.
\newblock {\em IEEE Trans. Parallel Distrib. Syst.}, 12(10):1094--1104.

\bibitem[\protect\astroncite{Mukherjee et~al.}{2018a}]{MBL17}
Mukherjee, D., Borst, S.~C., and {Van Leeuwaarden}, J. S.~H. (2018a).
\newblock {Asymptotically optimal load balancing topologies}.
\newblock {\em Proc. ACM Meas. Anal. Comput. Syst.}, 2(1):1--29.

\bibitem[\protect\astroncite{Mukherjee et~al.}{2018b}]{MBLW16-3}
Mukherjee, D., Borst, S.~C., van Leeuwaarden, J. S.~H., and Whiting, P.~A.
  (2018b).
\newblock {Universality of power-of-d load balancing in many-server systems}.
\newblock {\em Stoch. Syst.}, 8(4):265--292.

\bibitem[\protect\astroncite{Pang et~al.}{2007}]{PTRW07}
Pang, G., Talreja, R., and Whitt, W. (2007).
\newblock {Martingale proofs of many-server heavy-traffic limits for Markovian
  queues}.
\newblock {\em Probab. Surveys}, 4:193--267.

\bibitem[\protect\astroncite{Reiss et~al.}{2012}]{RTGKK12}
Reiss, C., Tumanov, A., Ganger, G.~R., Katz, R.~H., and Kozuch, M.~A. (2012).
\newblock Heterogeneity and dynamicity of clouds at scale: Google trace
  analysis.
\newblock In {\em Proceedings of the Third ACM Symposium on Cloud Computing},
  SoCC '12, pages 1--13. Association for Computing Machinery.

\bibitem[\protect\astroncite{Roughgarden and Tardos}{2002}]{RT02}
Roughgarden, T. and Tardos, E. (2002).
\newblock {How bad is selfish routing?}
\newblock {\em J. ACM}, 49(2):236--259.

\bibitem[\protect\astroncite{Stolyar}{1995}]{Stolyar95}
Stolyar, A.~L. (1995).
\newblock {On the stability of multiclass queueing networks: a relaxed
  sufficient condition via limiting fluid processes}.
\newblock {\em Markov Processes Relat.}, 1(4):491--512.

\bibitem[\protect\astroncite{Stolyar}{2005}]{Stolyar05}
Stolyar, A.~L. (2005).
\newblock {Optimal routing in output-queued flexible server systems}.
\newblock {\em Probab. Eng. Inf. Sci.}, 19(2):141--189.

\bibitem[\protect\astroncite{Stolyar}{2017}]{Stolyar17}
Stolyar, A.~L. (2017).
\newblock {Pull-based load distribution among heterogeneous parallel servers:
  the case of multiple routers}.
\newblock {\em Queueing Syst.}, 85(1):31--65.

\bibitem[\protect\astroncite{Tang and Subramanian}{2019a}]{TS19a}
Tang, D. and Subramanian, V.~G. (2019a).
\newblock Derandomized load balancing using random walks on expander graphs.

\bibitem[\protect\astroncite{Tang and Subramanian}{2019b}]{TS19}
Tang, D. and Subramanian, V.~G. (2019b).
\newblock {Random walk based sampling for load balancing in multi-server
  systems}.
\newblock {\em Proc. ACM Meas. Anal. Comput. Syst.}, 3(1):14:1--14:44.

\bibitem[\protect\astroncite{Tsitsiklis and Xu}{2013}]{TX13}
Tsitsiklis, J.~N. and Xu, K. (2013).
\newblock {Queueing system topologies with limited flexibility}.
\newblock In {\em Proc. SIGMETRICS '13}, pages 167--178.

\bibitem[\protect\astroncite{Tsitsiklis and Xu}{2017}]{TX17}
Tsitsiklis, J.~N. and Xu, K. (2017).
\newblock {Flexible queueing architectures}.
\newblock {\em Oper. Res.}, 65(5):1398--1413.

\bibitem[\protect\astroncite{Turner}{1998}]{Turner98}
Turner, S.~R. (1998).
\newblock {The effect of increasing routing choice on resource pooling}.
\newblock {\em Probab. Eng. Inf. Sci.}, 12(01):109--124.

\bibitem[\protect\astroncite{{van der Boor} et~al.}{2017}]{BBL17a}
{van der Boor}, M., Borst, S.~C., and van Leeuwaarden, J. S.~H. (2017).
\newblock {Load balancing in large-scale systems with multiple dispatchers}.
\newblock In {\em Proc. INFOCOM '17}.

\bibitem[\protect\astroncite{van~der Boor et~al.}{2017}]{BBLM18}
van~der Boor, M., Borst, S.~C., van Leeuwaarden, J. S.~H., and Mukherjee, D.
  (2017).
\newblock Scalable load balancing in networked systems: universality properties
  and stochastic coupling methods.

\bibitem[\protect\astroncite{{Van der Hofstad}}{2017}]{remco-book-1}
{Van der Hofstad}, R. (2017).
\newblock {\em {Random Graphs and Complex Networks}}, volume~1.
\newblock Cambridge University Press, Cambridge.

\bibitem[\protect\astroncite{Vargaftik et~al.}{2020}]{VKO20}
Vargaftik, S., Keslassy, I., and Orda, A. (2020).
\newblock {LSQ: Load balancing in large-scale heterogeneous systems with
  multiple dispatchers}.
\newblock {\em IEEE/ACM Trans. Netw.}, 28(3):1186--1198.

\bibitem[\protect\astroncite{Vvedenskaya et~al.}{1996}]{VDK96}
Vvedenskaya, N.~D., Dobrushin, R.~L., and Karpelevich, F.~I. (1996).
\newblock {Queueing system with selection of the shortest of two queues: An
  asymptotic approach}.
\newblock {\em Problemy Peredachi Informatsii}, 32(1):20--34.

\bibitem[\protect\astroncite{Wang et~al.}{2018a}]{WMSY18}
Wang, W., Maguluri, S.~T., Srikant, R., and Ying, L. (2018a).
\newblock Heavy-traffic delay insensitivity in connection-level models of data
  transfer with proportionally fair bandwidth sharing.
\newblock {\em SIGMETRICS Perform. Eval. Rev.}, 45(3):232--245.

\bibitem[\protect\astroncite{Wang et~al.}{2018b}]{WMSY18a}
Wang, W., Maguluri, S.~T., Srikant, R., and Ying, L. (2018b).
\newblock Heavy-traffic insensitive bounds for weighted proportionally fair
  bandwidth sharing policies.

\bibitem[\protect\astroncite{Wang et~al.}{2016}]{WZYTZ16}
Wang, W., Zhu, K., Ying, L., Tan, J., and Zhang, L. (2016).
\newblock {MapTask scheduling in MapReduce with data locality: Throughput and
  heavy-traffic optimality}.
\newblock {\em IEEE/ACM Trans. Netw.}, 24(1):190--203.

\bibitem[\protect\astroncite{Weng et~al.}{2020}]{WZS20}
Weng, W., Zhou, X., and Srikant, R. (2020).
\newblock Optimal load balancing in bipartite graphs.

\bibitem[\protect\astroncite{Xie et~al.}{2016}]{XYL16}
Xie, Q., Yekkehkhany, A., and Lu, Y. (2016).
\newblock {Scheduling with multi-level data locality: Throughput and
  heavy-traffic optimality}.
\newblock In {\em Proc. INFOCOM '16}, pages 1--9.

\bibitem[\protect\astroncite{Yekkehkhany et~al.}{2018}]{yekkehkhany2018gb}
Yekkehkhany, A., Hojjati, A., and Hajiesmaili, M.~H. (2018).
\newblock Gb-pandas: Throughput and heavy-traffic optimality analysis for
  affinity scheduling.
\newblock {\em SIGMETRICS Perform. Eval. Rev.}, 45(3):2--14.

\bibitem[\protect\astroncite{Yekkehkhany and Nagi}{2020}]{yekkehkhany2020blind}
Yekkehkhany, A. and Nagi, R. (2020).
\newblock Blind gb-pandas: A blind throughput-optimal load balancing algorithm
  for affinity scheduling.
\newblock {\em IEEE/ACM Transactions on Networking}, 28(3):1199--1212.

\bibitem[\protect\astroncite{Zhou et~al.}{2020}]{ZSW20}
Zhou, X., Shroff, N., and Wierman, A. (2020).
\newblock Asymptotically optimal load balancing in large-scale heterogeneous
  systems with multiple dispatchers.

\end{thebibliography}

\appendix
\section{Hard constraint on server-degrees}
\label{app:drg}

To establish Theorem~\ref{th:drg}, it suffices to prove that the sequence of random graphs as described in Theorem~\ref{th:drg} satisfies the conditions of proportional sparsity and subcriticality, almost surely.
Throughout this section, $\{G_N\}_{N \geq 1}$ will denote the sequence of random graphs as described in Theorem~\ref{th:drg}, where each server has degree $c(N)$. \\

\noindent
\textbf{Verification of proportional sparsity.}
We start by verifying the proportional sparsity condition, as stated in Proposition~\ref{prop:dregprop}. 
Lemma~\ref{lemma:chernoff_binom}--\ref{lemma:regulardispatchers} below provide necessary technical results for the proof of Proposition~\ref{prop:dregprop}.

\begin{lemma}
\label{lemma:chernoff_binom}
Let $X_N$ be a $\text{Binomial}(N, p(N))$ random variable. If $p(N) \to 0$ as $N \to \infty$, then for any fixed $\delta\in (0, 1/2)$, there exists $N_0 \geq 1$ and $a > 0$, such that for all $N \geq N_0$,
\begin{equation}
    \pro(X_N \geq \delta N) \leq \big( a (p(N))^\delta \big)^{N}.
\end{equation}
\end{lemma}

\begin{proof}
Fix $\delta\in (0, 1/2)$. Let $N_0 \geq 1$ such that $p(N) \leq \delta$ for all $N \geq N_0$. 
By the Chernoff bound for binomials \cite[Equation 2.4]{JLR00}, we get for all $N \geq N_0$,
\begin{equation}
\begin{split}
    \pro(&X_N \geq \delta N) \\
    &\leq \exp\Big( -N \Big( (\delta + p(N)) \ln\Big(\frac{\delta + p(N)}{p(N)}\Big) 
    - (1 - p(N) - \delta) \ln\Big(\frac{1 - p(N)}{1 - p(N) - \delta}\Big) \Big) \Big) \\
    &\leq \exp\Big( -N \Big( \delta \ln\Big(\frac{\delta}{p(N)}\Big) - (1 - \delta) \ln\Big(\frac{1 - p(N)}{1 - 2 \delta}\Big) \Big) \Big) \\
    &= \exp\Big( -\delta \ln\Big(\frac{1}{p(N)}\Big) + \Big( \delta \ln\Big(\frac{1}{\delta}\Big) + (1 - \delta) \ln\Big(\frac{1}{1 - 2 \delta}\Big) \Big) \Big)^{N} 
    = \big( a (p(N))^\delta \big)^{N},
\end{split}
\end{equation}
where $a = \exp\left( \delta \ln\left(\frac{1}{\delta}\right) + (1 - \delta) \ln\left(\frac{1}{1 - 2 \delta}\right) \right)$.
\end{proof}
Define the neighborhood of a server $v \in V_N$ as $\cN_v := \{ w \in W_N \mid (v, w) \in E_N \}$. The next lemma provides sufficient conditions on the server-degree $c(N)$ to ensure that the number of edges from the servers to any subset of dispatchers is close to its mean, almost surely.

\begin{lemma}
\label{lemma:regularservers}
For any $X \subseteq W_N$ and $\varepsilon>0$, define
\begin{equation}
B_N(X) := \left\{ v \in V_N \quad\mid\quad \left\lvert \lvert \cN_v \cap X \rvert - \frac{c(N) \lvert X \rvert}{M(N)} \right\rvert \geq \frac{\varepsilon c(N) \lvert X \rvert}{M(N)} \right\}.
\end{equation}
If $c(N) = \omega(1)$ and $c(N) = \omega\left(\frac{M(N)}{N}\right)$, then for all $0 < \varepsilon \leq 3 / 2$, $0 < \delta < 1/2$, and $\eta > 0$, almost surely there exists $N_0 \geq 1$ such that, for all $N \geq N_0$,
$
    \lvert B_N(X) \rvert \leq \delta N,
$
for all $X \subseteq W_N$ with $\lvert X \rvert \geq \eta M(N)$.
\end{lemma}

\begin{proof}
Fix $0 < \varepsilon \leq 3 / 2$, $0 < \delta < 1/2$ and $\eta > 0$. The number of edges between $v \in V_N$ and $X \subseteq W_N$ is a Hypergeometric random variable. By the Chernoff bound for Hypergeometric random variables \cite[Corollary 2.3]{JLR00}, we have for any $X \subseteq W_N$ with $\lvert X \rvert \geq \eta M(N)$,
\begin{equation}
\begin{gathered}
    \pro\left( \left\lvert \lvert \cN_v \cap X \rvert - \frac{c(N) \lvert X \rvert}{M(N)} \right\rvert \geq \frac{\varepsilon c(N) \lvert X \rvert}{M(N)} \right)
    \leq 2 \exp\left( -\frac{\varepsilon^2 c(N) \lvert X \rvert}{3 M(N)} \right) \\
    \leq 2 \exp\left( -\frac{\varepsilon^2 \eta c(N)}{3} \right) =: p(N).
\end{gathered}
\end{equation}
For any fixed $X$, the events in the definition of $B_N(X)$ are independent over $v \in V_N$. The random variable $\lvert B_N(X) \rvert$ is therefore stochastically upper bounded by a Binomial $(N, p(N))$ random variable for all $X$. By Lemma~\ref{lemma:chernoff_binom} and because $c(N) \to \infty$ as $N \to \infty$, there exists $N_1 \geq 1$ and $a > 0$ such that,
\begin{equation}
\begin{split}
    \pro\big(\lvert B_N(X) \rvert \geq \delta N\big)
    \leq \left( a p(N)^\delta \right)^N 
    &= \exp\Big( - N \Big(\varepsilon^2 \eta \delta c(N)/3 - \ln(2a) \Big) \Big)\\
    &\leq \exp\big( -\varepsilon^2 \eta \delta c(N) N/6 \big),
\end{split}
\end{equation}
for all $N \geq N_1$ and $X \subseteq W_N$. Because $c(N) N / M(N) \to \infty$ and $c(N) \to \infty$ as $N \to \infty$, there exists $N_2 \geq 1$ such that,
\begin{equation}
\begin{gathered}
    \pro(\exists\ X \subseteq W_N\text{ with } \lvert X \rvert \geq \eta M(N)\text{ such that } \lvert B_N(X) \rvert \geq \delta N)
    \leq 2^{M(N)} \exp\left( -\frac{\varepsilon^2 \eta \delta c(N) N}{6} \right) \\
    = \exp\left( \ln(2) M(N) - \frac{\varepsilon^2 \eta \delta c(N) N}{6} \right)
    \leq \exp\left( -\frac{\varepsilon^2 \eta \delta c(N) N}{12} \right)
    \leq \exp(-N)
\end{gathered}
\end{equation}
for all $N \geq N_2$. Hence, the proof is completed by the first Borel-Cantelli lemma.
\end{proof}

We will now use Lemma~\ref{lemma:regularservers} to prove that for \emph{most dispatchers}, the number of servers sharing an edge with the dispatcher is close to the mean.

\begin{lemma}
\label{lemma:regulardispatchers}
Define
\begin{equation}
A_N := \big\{ w \in W_N \mid \lvert \cN_w \rvert \leq (1 - \varepsilon) \mathbb{E}[ \lvert \cN_w \rvert ] \big\}.
\end{equation}
If $c(N) = \omega(1)$ and $c(N) = \omega\left(\frac{M(N)}{N}\right)$, then for all $\varepsilon > 0$ and $\delta > 0$, almost surely there exists $N_0 \geq 1$ such that,
$
     \lvert A_N \rvert \leq \delta M(N),
$
for all $N \geq N_0$.
\end{lemma}

\begin{proof}
Fix $\varepsilon > 0$ and $\delta > 0$. Assume for the sake of contradiction that there is a sequence in $\mathbb{N}$ for which $\lvert A_N \rvert \geq \delta M(N)$ for all $N$ in this sequence. The number of edges between $V_N$ and $A_N$ is then
\begin{equation}
\label{eq:maxedges}
    \sum_{w \in A_N} \lvert \cN_w \rvert
    \leq \frac{(1 - \varepsilon) c(N) N \lvert A_N \rvert}{M(N)}.
\end{equation}
However, by Lemma \ref{lemma:regularservers}, there exists $N_0 \geq 1$ such that,
\begin{equation}
\begin{split}
    &\sum_{v \in V_N} \lvert \cN_v \cap A_N \rvert
    \geq \sum_{v \in B_N(A_N)^c} \lvert \cN_v \cap A_N \rvert \\
    &\geq (N - \lvert B_N(A_N) \rvert) \frac{(1 - \varepsilon_1) c(N) \lvert A_N \rvert}{M(N)}
    \geq (1 - \delta_1) \frac{(1 - \varepsilon_1) c(N) N \lvert A_N \rvert}{M(N)},
\end{split}
\end{equation}
for all $N \geq N_0$ which contradicts equation \eqref{eq:maxedges} for $\varepsilon_1$ and $\delta_1$ small enough.
\end{proof}


\begin{proposition}
\label{prop:dregprop}
Assume $c(N) = \omega(1)$ and $c(N) = \omega\left(\frac{M(N)}{N}\right)$. 
Then the sequence of graphs $\{G_N\}_{N\geq 1}$ is proportionally sparse almost surely.
\end{proposition}

\begin{proof}
Let $\varepsilon_1, \varepsilon_2 > 0$ and define 
\begin{align*}
    A_N &:= \big\{ w \in W_N \;\mid\; \lvert \cN_w \rvert \leq (1 - \varepsilon_1) \mathbb{E}[ \lvert \cN_w \rvert ] \big\}\\
    B_N(X) &:= \left\{ v \in V_N \quad\mid\quad \left\lvert \lvert \cN_v \cap X \rvert - \frac{c(N) \lvert X \rvert}{M(N)} \right\rvert \geq \frac{\varepsilon_2 c(N) \lvert X \rvert}{M(N)} \right\}\text{ for } X \subseteq W_N.
\end{align*}
Fix $\varepsilon > 0$ and $\delta > 0$. The goal will be to prove proportional sparsity, that is, there exists $N_0 \geq 1$ such that
\begin{equation}
    \sup_{U \subseteq V_N} \left\lvert \left\{ w \in W_N \mid \left\lvert \frac{\lvert \cN_w \cap U \rvert}{\lvert \cN_w \rvert} - \frac{\lvert U \rvert}{N} \right\rvert \geq \varepsilon \right\} \right\rvert \leq 2 \delta M(N) \text{ for } N \geq N_0.
\end{equation}
Lemma \ref{lemma:regulardispatchers} showed that $\lvert A_N \rvert \leq \delta M(N)$ for $N$ large enough. It is therefore sufficient to prove that there exists $N_1 \geq 1$ such that
\begin{equation}
    \sup_{U \subseteq V_N} \left\lvert \left\{ w \in A_N^c \mid \left\lvert \frac{\lvert \cN_w \cap U \rvert}{\lvert \cN_w \rvert} - \frac{\lvert U \rvert}{N} \right\rvert \geq \varepsilon \right\} \right\rvert \leq \delta M(N) \text{ for } N \geq N_1.
\end{equation}
Consider a subset $U \subseteq V_N$. We proceed by contradiction. First, consider the case that there is a subset of dispatchers which is over connected to $U$. In other words, assume that there exists a subset $X \subseteq A_N^c \subseteq W_N$ with $\lvert X \rvert \geq \delta M(N)$ such that for all $w \in X$,
\begin{equation}
    \frac{\lvert \cN_w \cap U \rvert}{\lvert \cN_w \rvert} \geq \frac{\lvert U \rvert}{N} + \varepsilon.
\end{equation}
We now distinguish two cases based on the cardinality of $U$. Choose $\eta := \varepsilon \delta / 2$ and assume that $\lvert U \rvert \geq \eta N$. By Lemma \ref{lemma:regularservers}, the number of edges between $U$ and $X$ is
\begin{equation}
\begin{gathered}
\label{eq:minimumedges_1}
    \sum_{v \in U \cap B_N(X)} \lvert \cN_v \cap X \rvert + \sum_{v \in U \cap B_N(X)^c} \lvert \cN_v \cap X \rvert
    \leq c(N) \lvert B_N(X) \rvert + \frac{(1 + \varepsilon_2) c(N) \lvert X \rvert \lvert U \rvert}{M(N)} \\
    \leq \frac{(1 + 2 \varepsilon_2) c(N) \lvert X \rvert \lvert U \rvert}{M(N)},
\end{gathered}
\end{equation}
for all $N$ large enough such that $\lvert B_N(X) \rvert \leq \varepsilon_2 \delta \eta N$. At the same time, the number of edges between $U$ and $X$ is
\begin{equation}
\begin{gathered}
    \sum_{w \in X} \lvert \cN_w \cap U \rvert
    \geq \sum_{w \in X} \lvert \cN_w \rvert \left( \frac{\lvert U \rvert}{N} + \varepsilon \right)
    \geq (1 - \varepsilon_1) \frac{c(N) N}{M(N)} \left( \frac{\lvert U \rvert}{N} + \varepsilon \right) \lvert X \rvert \\
    = (1 - \varepsilon_1) (1 + \varepsilon) \frac{c(N) \lvert X \rvert \lvert U \rvert}{M(N)}
\end{gathered}
\end{equation}
which contradicts equation \eqref{eq:minimumedges_1} for $\varepsilon_1$ and $\varepsilon_2$ small enough.

Assume now that $\lvert U \rvert < \eta N$. The number of edges between $U$ and $X$ is
\begin{equation}
\label{eq:minimumedges_2}
    \sum_{v \in U} \lvert \cN_v \cap X \rvert \leq \eta c(N) N = \varepsilon \delta c(N) N / 2.
\end{equation}
Then, the number of edges between $U$ and $X$ is
\begin{equation}
\begin{gathered}
    \sum_{w \in X} \lvert \cN_w \cap U \rvert
    \geq \sum_{w \in X} \lvert \cN_w \rvert \left( \frac{\lvert U \rvert}{N} + \varepsilon \right)
    \geq (1 - \varepsilon_1) \frac{c(N) N}{M(N)} \left( \frac{\lvert U \rvert}{N} + \varepsilon \right) \lvert X \rvert \\
    \geq \varepsilon \delta (1 - \varepsilon_1) c(N) N 
\end{gathered}
\end{equation}
which contradicts equation \eqref{eq:minimumedges_2} for $\varepsilon_1$ small enough.

Second, consider the case that there is a subset of dispatchers which is under connected to $U$. In other words, assume that there exists $X \subseteq A_N^c \subseteq W_N$ with $\lvert X \rvert \geq \delta M(N)$ such that for all $w \in X$,
\begin{equation}
    \frac{\lvert \cN_w \cap U \rvert}{\lvert \cN_w \rvert}
    \leq \frac{\lvert U \rvert}{N} - \varepsilon.
\end{equation}
Then we consider the complement of $U$ to find
\begin{equation}
    \frac{\lvert \cN_w \cap U^c \rvert}{\lvert \cN_w \rvert}
    = 1 - \frac{\lvert \cN_w \cap U \rvert}{\lvert \cN_w \rvert}
    \geq 1 - \frac{\lvert U \rvert}{N} + \varepsilon
    = \frac{\lvert U^c \rvert}{N} + \varepsilon,
\end{equation}
which leads to a contradiction in the same way as before. Hence, the sequence $\{G_N\}_{N\geq 1}$ as described in Theorem~\ref{th:drg} is proportionally sparse almost surely.
\end{proof}

\noindent
\textbf{Verification of subcriticality.}
We now prove that $\{G_N\}_{N\geq 1}$ satisfies the subcriticality condition in Definition~\ref{def:subcritical}, almost surely. 
Note that it is enough to verify that $\{G_N\}_{N\geq 1}$ satisfies~\eqref{eq:subcritical} for some choice of $\gamma_w^U(v)$'s.
In particular, if $\gamma_w^U(v)$ is the uniform distribution, i.e. $\gamma_w^U(v) = 1/\min(d, \lvert \cN_w \rvert)$ for $v \in U$, then the condition reduces to
\begin{equation}
\label{eq:subcrituniform}
    \limsup_{N \to \infty} \max_{v \in V_N} \frac{N}{M(N)} \sum_{w \in W_N} \frac{\mathbbm{1}\{ (v, w) \in E_N \}}{\lvert \cN_w \rvert} \leq 1.
\end{equation}
In the rest of this section, we will verify that $\{G_N\}_{N\geq 1}$ satisfies \eqref{eq:subcrituniform} almost surely.

\begin{remark}\normalfont
The condition of subcriticality as in~\eqref{eq:subcrituniform} is similar to Condition 1(ii) in \cite{BMW17}. Translated to the current framework, Condition 1(ii) in \cite{BMW17} reads
\begin{equation}
    \max_{v \in V_N} \frac{N}{M(N)} \left\lvert \sum_{w \in W_N} \frac{\mathbbm{1}\{ (v, w) \in E_N \}}{\lvert \cN_w \rvert} - 1 \right\rvert \to 0 \text{ as } N \to \infty.
\end{equation}
The subcriticality condition in~\eqref{eq:subcrituniform} is less restrictive than the condition in \cite{BMW17}.
\end{remark}


\begin{lemma}
\label{lemma:drgmax}
Assume $c(N) = \omega(M(N) \ln(N) / N)$.  Then the graph sequence $\{G_N\}_{N\geq 1}$ satisfies~\eqref{eq:subcrituniform}, and hence the subcriticality condition, almost surely.
\end{lemma}

\begin{proof}
Let $\zeta(N) := \sqrt{\ln(N) c(N) N / M(N)}$ and notice that 
\begin{equation}
\frac{\zeta(N)}{c(N) N / M(N)} = \sqrt{\frac{\ln(N)}{c(N) N / M(N)}} \to 0 \text{ as } N \to \infty.
\end{equation}

We know that $\lvert \cN_w \rvert$ is Binomial with parameters $N$ and $c(N) / M(N)$. By the Chernoff bound for Binomials~\cite[Theorem 2.1]{JLR00},
\begin{equation}
\begin{split}
    \mathbb{P}\left( \lvert \cN_w \rvert \leq \frac{c(N) N}{M(N)} - a \zeta(N) \right) 
    &\leq \exp\left( -\frac{a^2 \zeta(N)^2}{2 c(N) N / M(N)} \right) \\
    &\leq \exp\left( -\frac{a^2 \ln(N)}{2} \right) = N^{ -a^2 / 2 }.
\end{split}
\end{equation}
Hence,
\begin{equation}
    \mathbb{P} \left( \exists w \in W_N \text{ such that } \lvert \cN_w \rvert \geq \frac{c(N) N}{M(N)} - a \zeta(N) \right) \leq N^{1 - a^2 / 2}
\end{equation}
As the probabilities are summable over $N$ for arbitrary $a^2 > 4$, the first Borel-Cantelli lemma proves the events can only happen finitely many times. Hence, there exists $N_0 \geq 1$ such that
\begin{equation}
    \lvert \cN_w \rvert \geq c(N) N / M(N) - a \zeta(N) \text{ for all } w \in W_N \text{ and } N \geq N_0. 
\end{equation}
Therefore,
\begin{align*}
    &\max_{v \in V_N} \frac{N}{M(N)} \sum_{w \in W_N} \frac{\mathbbm{1}\{ (v, w) \in E_N \}}{\lvert \cN_w \rvert}
    \leq \frac{N}{M(N)} \frac{c(N)}{\min_{W \in V_N} \lvert \cN_w \rvert} \\
    &\leq \frac{N}{M(N)} \frac{c(N)}{c(N) N / M(N) - a \zeta(N)}
    = 1 + \frac{a \zeta(N)}{c(N) M(N) / N - a \zeta(N)} \to 1 \text{ as } N \rightarrow \infty.
\end{align*}
\end{proof}

\begin{proof}[Proof of Theorem~\ref{th:drg}]
The proof of Theorem~\ref{th:drg} follows immediately from Proposition~\ref{prop:dregprop} and Lemma~\ref{lemma:drgmax}.
\end{proof}

\section{Inhomogeneous levels of flexibility}
\label{app:erg}

Throughout this section, $\{G_N\}_{N\geq 1}$ will denote the sequence of random graphs as described in Theorem~\ref{th:erg}, where $(p_w(N))_{w\in W_N}$ is the connection probability vector. 
To establish Theorem~\ref{th:erg}, it suffices to prove that $\{G_N\}_{N\geq 1}$ satisfies the conditions of proportional sparsity and subcriticality, almost surely.
\\

\noindent
\textbf{Verification of proportional sparsity.}
We start by verifying the proportionally sparsity condition.  
Define $\bar{p}(N) := \min_{w \in W_N} p_w(N)$. 


\begin{lemma}
\label{lemma:ergprop}
Assume $\bar{p}(N) = \omega(1 / N)$ and $\bar{p}(N) = \omega(1 / M(N))$. Then the sequence of graphs $\{ G_N \}_{N \geq 1}$ is proportionally sparse almost surely.
\end{lemma}

\begin{proof}
Fix $\varepsilon > 0$, $w \in W_N$ and $U \subseteq V_N$. Let $B_w(U)$ denote the event that a dispatcher $w \in W_N$ is bad with respect to the set $U \subseteq V_N$, defined as
\begin{equation}
    B_w(U) := \left\{ \left\lvert \frac{\lvert \cN_w \cap U \rvert}{\lvert \cN_w \rvert} - \frac{\lvert U \rvert}{N} \right\rvert \geq \varepsilon \right\}.
\end{equation}
Define $\zeta_w(N) := p_w(N) N$ and $\eta := \lvert U \rvert / N$. By the law of total probability,
\begin{equation}
\begin{gathered}
\label{eq:erg_totalsplit}
    \mathbb{P}\left( B_w(U) \right) \\
    \leq \mathbb{P}\big( B_w(U) \;\mid\; \big\lvert \lvert \cN_w \cap U \rvert - \eta \zeta_w(N) \big\rvert < \varepsilon_1 \zeta_w(N) \text{ and } \big\lvert \lvert \cN_w \rvert - \zeta_w(N) \big\rvert < \varepsilon_2  \zeta_w(N) \big) \\
    + \mathbb{P}\left( \big\lvert \lvert \cN_w \cap U \rvert - \eta \zeta_w(N) \big\rvert \geq \varepsilon_1 \zeta_w(N)) + \mathbb{P}(\big\lvert \lvert \cN_w \rvert - \zeta_w(N) \big\rvert \geq \varepsilon_2 \zeta_w(N) \right).
\end{gathered}
\end{equation}
We will bound each terms on the right-hand-side of \eqref{eq:erg_totalsplit}. The first term becomes equal to zero by choosing $\varepsilon_1 = \varepsilon / 6$ and $\varepsilon_2 = \min(\varepsilon / 6, 1/2)$ such that
\begin{equation}
    \frac{\lvert \cN_w \cap U \rvert}{\lvert \cN_w \rvert} - \eta \leq \frac{\eta \zeta_w(N) + \varepsilon_1 \zeta_w(N)}{\zeta_w(N) - \varepsilon_2 \zeta_w(N)} - \eta \leq \frac{\varepsilon_1 + \varepsilon_2}{1 - \varepsilon_2} < \varepsilon
\end{equation}
and
\begin{equation}
    \frac{\lvert \cN_w \cap U \rvert}{\lvert \cN_w \rvert} - \eta \geq \frac{\eta \zeta_w(N) - \varepsilon_1 \zeta_w(N)}{\zeta_w(N) + \varepsilon_2 \zeta_w(N)} - \eta \geq -\frac{\varepsilon_1 + \varepsilon_2}{1 + \varepsilon_2} > -\varepsilon.
\end{equation}
By applying the Chernoff bound for Binomials \cite[Theorem 2.1]{JLR00} on the second term we obtain
\begin{equation}
\begin{split}
    &\mathbb{P}\left( \lvert \lvert \cN_w \cap U \rvert - \eta \zeta_w(N) \rvert \geq (\varepsilon_1 / \eta) \eta \zeta_w(N) \right)\\
    &\leq 2 \exp\left( -\eta \zeta_w(N) (\varepsilon_1 / \eta)^2/3 \right)
    \leq 2 \exp\left( -\bar{p}(N) N \varepsilon_1^2/3 \right),
\end{split}
\end{equation}
and applying the Chernoff bound \cite[Theorem 2.1]{JLR00} on the third term yields
\begin{equation}
    \mathbb{P}\left( \lvert \lvert \cN_w \rvert - \zeta_w(N) \rvert \geq \varepsilon_2  p_w(N) N \right)
    \leq 2 \exp\left( -\frac{\zeta_w(N) \varepsilon_2^2}{3} \right)
    \leq 2 \exp\left( -\frac{\bar{p}(N) N \varepsilon_2^2}{3} \right).
\end{equation}
Therefore,
\begin{equation}
\label{eq:bupper}
    \mathbb{P}(B_w(U)) \leq 2 \exp\left( -\frac{\bar{p}(N) N \varepsilon_1^2}{3} \right) + 2 \exp\left( -\frac{\bar{p}(N) N \varepsilon_2^2}{3} \right),
\end{equation}
uniformly over $w \in W_N$ and $U \subseteq V_N$. Fix $0 < \delta < 1/2$ and define $\alpha(\delta)$ and $\beta(\delta)$ as
\begin{align}
    \alpha(\delta) := \frac{4}{\delta}, &&
    \beta(\delta) := \exp\left( -\frac{2}{\delta} \left (\delta \ln\left(\frac{1}{\delta}\right) + (1 - \delta) \ln\left(\frac{1}{1 - 2\delta}\right) \right) \right).
\end{align}
Because $\bar{p}(N) N \to \infty$ as $N \to \infty$, equation \eqref{eq:bupper} converges to zero as $N$ becomes large. Hence, for all $\delta > 0$ there exists $N_0 \geq 1$ such that $\mathbb{P}(B_w(U)) \leq \beta(\delta)$ for all $N \geq N_0$. Similarly, because $\bar{p}(N) M(N) \to \infty$ as $N \to \infty$, it follows that for all $\delta > 0$ there exists $N_1 \geq 1$ such that $\mathbb{P}(B_w(U)) \leq \exp\left( -\frac{\alpha(\delta) N}{M(N)} \right)$ for all $N \geq N_1$. As the events $\{ B_w(U) \}_{w \in W_N}$ are independent, for all $\delta > 0$ and $N \geq \max(N_0, N_1)$, the sum of their indicators $\sum_{w \in W_N} \mathbbm{1}(B_w(U))$ is stochastically dominated by a Binomial $(M(N), \theta)$ random variable, where
\begin{equation}
\theta := \min\left( \exp\left( - \frac{\alpha(\delta) N}{M(N)} \right), \beta(\delta) \right).
\end{equation}
By the Chernoff bound for Binomials  \cite[Theorem 2.1]{JLR00},
\begin{equation}
\begin{split}
    &\mathbb{P}\left( \sum_{w \in W_N} \mathbbm{1}(B_w(U)) \geq \delta M(N) \right) \\
   & \leq \exp\left( -M(N) \left( \delta \ln\left(\frac{\delta}{\theta}\right) - (1 - \delta) \ln\left(\frac{1 - \theta}{1 - 2\delta}\right) \right) \right) \\
    &\leq \exp\left( -\frac{\delta M(N)}{2} \ln\left(\frac{1}{\theta}\right) + M(N) \left( \delta \ln\left(\frac{1}{\delta}\right) + (1 - \delta) \ln\left(\frac{1}{1 - 2\delta}\right) - \frac{\delta}{2} \ln\left(\frac{1}{\theta}\right) \right) \right) \\
    &\leq \exp\left( - \frac{\delta M(N)}{2} \frac{\alpha(\delta) N}{M(N)} + M(N) \left( \delta \ln\left(\frac{1}{\delta}\right) + (1 - \delta) \ln\left(\frac{1}{1 - 2\delta}\right) - \frac{\delta}{2} \ln\left(\frac{1}{\beta(\delta)}\right) \right) \right) \\
    &\leq \exp\left( -2 N \right) \\
\end{split}
\end{equation}
for all $N \geq \max(N_0, N_1)$, and hence,
\begin{equation}
\begin{split}
   & \mathbb{P}\left( \sup_{U \subseteq V_N} \sum_{w \in W_N} \mathbbm{1}(B_w(U)) \geq \delta M(N) \right) \\
    &\qquad\leq
    \sum_{U \subseteq V_N} \mathbb{P}\left( \sum_{w \in W_N} \mathbbm{1}(B_w(U)) \geq \delta M(N) \right) 
    \leq \exp (-N).
\end{split}
\end{equation}
As the probabilities are summable over $N$, the first Borel-Cantelli lemma proves that the sequence of random graphs is proportionally sparse almost surely.
\end{proof}

\ \\
\noindent
\textbf{Verification of subcriticality.}
We now prove that $\{G_N\}_{N\geq 1}$ satisfies the subcriticality condition in Definition~\ref{def:subcritical}, almost surely. 
Note, as in Appendix~\ref{app:drg}, we will verify that $\{G_N\}_{N\geq 1}$ satisfies \eqref{eq:subcrituniform} almost surely. Recall that $\bar{p}(N) = \min_{w \in W_N} p_w(N)$.

\begin{lemma}
\label{lemma:ergmax}
Assume 
\begin{align*}
    \bar{p}(N) &= \omega((\ln(M(N)) + \ln(N)) / N),\\
    \lVert \boldsymbol{p}(N)^{-1} \rVert_2^{-2} &= \omega(\ln(N) / M(N)^2).
\end{align*}
 Then the sequence of graphs $\{ G_N \}_{N \geq 1}$ satisfies the subcriticality condition almost surely.
\end{lemma}

\begin{proof}
Let $\zeta_w(N) := \sqrt{(\ln(M(N)) + \ln(N)) p_w(N) N}$ and notice that
\begin{equation}
\frac{\zeta_w(N)}{p_w(N) N} \leq \sqrt{\frac{\ln(M(N)) + \ln(N)}{\bar{p}(N) N}} \to 0 \text{ as } N \to \infty.
\end{equation}
We know that $\lvert \cN_w \rvert$ is a Binomial $(N, p_w(N))$ random variable. By the Chernoff bound for Binomials \cite[Theorem 2.1]{JLR00},
\begin{equation}
\begin{gathered}
    \mathbb{P}(\lvert \cN_w \rvert \geq p_w(N) N - c_1 \zeta_w(N))
    \leq \exp\left( -\frac{c_1^2 \zeta_w(N)^2}{2 p_w(N) N} \right) \\
    \leq \exp\left( -\frac{c_1^2 (\ln(M(N)) + \ln(N))}{2} \right) = M(N)^{-c_1^2 / 2} \cdot N^{-c_1^2 / 2}.
\end{gathered}
\end{equation}
Hence,
\begin{equation}
    \mathbb{P}\big(\exists w \in W_N \text{ such that } \lvert \cN_w \rvert \geq p_w(N) N - c_1 \zeta_w(N)\big) \leq M(N)^{1 - c_1^2 / 2} \cdot N^{ -c_1^2 / 2}
\end{equation}
As the probabilities are summable over $N$ for arbitrary $c_1^2 > 2$, due to the first Borel-Cantelli lemma, almost surely, there exists $N_0 \geq 1$ such that for all $N \geq N_0$,
\begin{equation}
    \lvert \cN_w \rvert \geq p_w(N) N - c_1 \zeta_w(N) \text{ for all } w \in W_N. 
\end{equation}
Define the functions $f_v$ and $g$ for $v \in V_N$ as
\begin{equation}
\begin{aligned}
    f_v(\boldsymbol{z}) = \frac{N}{M(N)} \sum_{w \in W_N} \frac{z_w}{\lvert \cN_w \setminus \{ v \} \rvert + z_w}, &&
    g(\boldsymbol{z}) = \frac{N}{M(N)} \sum_{w \in W_N} \frac{z_w}{p_w(N) N - c_1 \zeta_w(N)},
\end{aligned}
\end{equation}
for vectors $\boldsymbol{z} \in \{ 0, 1 \}^{M(N)}$. Note that for $N \geq N_0$,
\begin{equation}
\begin{gathered}
    f_v(\boldsymbol{z}) - g_v(\boldsymbol{z})
    = \frac{N}{M(N)} \sum_{w \in W_N} z_w \left( \frac{1}{\lvert \cN_w \setminus \{ v \} \rvert + 1} - \frac{1}{p_w(N) N - c_1 \zeta_w(N)} \right) \\
    \leq \frac{N}{M(N)} \sum_{w \in W_N} z_w \left( \frac{1}{\lvert \cN_w \rvert} - \frac{1}{p_w(N) N - c_1 \zeta_w(N)} \right)
    \leq 0,
\end{gathered}
\end{equation}
and therefore $f_v(\boldsymbol{z}) \leq g_v(\boldsymbol{z})$. Define the edge indicators $Z_{v,w} = \mathbbm{1}\{ (v, w) \in E_N \}$ for $v \in V_N$ and $w \in W_N$. Note that the random variables $\{ Z_{v,w} \}_{w \in W_N}$ are independent. Furthermore, there exists $N_1 \geq 1$ such that if two vectors $\boldsymbol{z}, \boldsymbol{z}' \in \{ 0, 1 \}^{M(N)}$ differ only in the $w$-th coordinate then
\begin{equation}
    \lvert g_v(\boldsymbol{z}) - g_v(\boldsymbol{z}') \rvert
    = \frac{N}{M(N)} \frac{1}{p_w(N) N - c_1 \zeta_w(N)}
    \leq \frac{2}{p_w(N) M(N)},
\end{equation}
for $N \geq N_1$. Finally, note that
\begin{equation}
    \mathbb{E}\left[ g_v\left( \{ Z_{v,w} \}_{w \in W_N} \right) \right] \to 1 \text{ as } N \to \infty.
\end{equation}
We apply Azuma-Hoeffding\cite[Corollary 2.27]{JLR00} and use the condition $\lVert \boldsymbol{p}(N)^{-1} \rVert_2^{-2} = \omega(\ln(N) / M(N)^2)$ to obtain
\begin{equation}
\begin{gathered}
    \mathbb{P}\left( f_v\left( \{ Z_{v,w} \}_{w \in W_N} \right) \geq 1 + 2 \varepsilon \right)
    \leq \mathbb{P}\left( g_v\left( \{ Z_{v,w} \}_{w \in W_N} \right) \geq \mathbb{E}\left[ g_v\left( \{ Z_{v,w} \}_{w \in W_N} \right) \right] + \varepsilon \right) \\
    \leq \exp\left( -\frac{\varepsilon^2}{2 \sum_{w \in W_N} \frac{4}{p_w(N)^2 M(N)^2}} \right)
    = \exp\left( -\frac{\varepsilon^2 M(N)^2}{8 \lVert \boldsymbol{p}(N)^{-1} \rVert_2^2} \right)
    \leq \exp\left( -\frac{c_2 \varepsilon^2 \ln(N)}{8} \right)
\end{gathered}
\end{equation}
for all $c_2 > 0$ and $N$ large enough. Therefore,
\begin{equation}
    \mathbb{P}\left( \exists v \in V_N \text{ such that } \max_{v \in V_N} \frac{N}{M(N)} \sum_{w \in W_N} \frac{\mathbbm{1}\{ (v, w) \in E_N \}}{\lvert \cN_w \rvert} \geq 1 + 2 \varepsilon \right)
    \leq N^{1 -\frac{c_2 \varepsilon^2}{8}}.
\end{equation}
As the probabilities are summable over $N$ for $c_2$ large enough, the first Borel-Cantelli lemma proves that the graph sequence satisfies the subcriticality condition almost surely.
\end{proof}

\begin{proof}[Proof of Theorem~\ref{th:erg}]
The proof of Theorem~\ref{th:erg} follows immediately from Lemmas~\ref{lemma:ergprop} and~\ref{lemma:ergmax}.
\end{proof}

\section{Global stability analysis}
\label{app:proofglobalstability}

The proof of global stability is based on the proof of Theorem 3.6 in \cite{Mitzenmacher01}. Recall that $\Psi_\omega(t) := \sum_{i = 1}^\infty \omega_i \lvert q_i^*(t) - q_i^*(\infty) \rvert$ and $i_0 := \min\{ i \geq 1 \mid \lambda (2 q_i^*(\infty) + 1) < \frac{1 + \lambda}{2} \}$.

\begin{proof}[Proof of \ref{th:globalstability}]
To establish the theorem, it is sufficient to show that there exists $1 < r < 2 / (1 + \lambda)$ such that $\boldsymbol{\omega} \in \mathbb{R}^\infty$ satisfies the condition from the proof of Theorem 3.6 in \cite{Mitzenmacher01}. That is, if there exists $\delta > 0$ such that
\begin{equation}
\label{eq:globalstabcon}
    \omega_{i + 1} \leq \omega_i + \frac{\omega_i (1 - \delta) - \omega_{i - 1}}{\lambda (2 q_i^*(\infty) + 1)} \quad\text{for}\quad i \geq 1,
\end{equation}
then $\Psi_\omega$ converges exponentially to zero if $\Psi(0) < \infty$. As suggested in the proof of Theorem 3.6 in \cite{Mitzenmacher01}, the weights are broken up into two subsequences starting with $\omega_0 = 0$ and $\omega_1 = 1$. For $1 \leq i \leq i_0 - 1$, we set
\begin{equation}
    \omega_{i + 1} := \omega_i + \frac{\omega_i (1 - \delta) - \omega_{i - 1}}{3} \leq \omega_i + \frac{\omega_i (1 - \delta) - \omega_{i - 1}}{\lambda (2 q_i^*(\infty) + 1)}.
\end{equation}
Note that the subsequence $\omega_0, \omega_1, \dots, \omega_{i_0}$ consists of finitely many terms. Hence, there exists $\delta_0 > 0$ small such that this subsequence is increasing for $\delta \leq \delta_0$. Applying equation \eqref{eq:globalstabcon} to the $(i_0 + 1)$-th term yields
\begin{equation}
    \omega_{i_0 + 1} := \omega_{i_0} r \leq \omega_{i_0} + \frac{\omega_{i_0} (1 - \delta) - \omega_{i_0 - 1}}{\lambda (2 \pi_{i_0} + 1)}.
\end{equation}
from which it follows that
\begin{equation}
    r \leq 1 + \frac{1}{\omega_{i_0}} \frac{\omega_{i_0} (1 - \delta) - \omega_{i_0 - 1}}{\lambda (2 \pi_{i_0} + 1)} =: R(\delta).
\end{equation}
Note that $R(\delta)$ increases as $\delta$ decreases. Hence, there exists $\delta_1 > 0$ small such that $R(\delta) > 1 + \epsilon$ for $\delta \leq \delta_1$ and some $\epsilon > 0$. Define
\begin{equation}
    r(\delta) := \frac{1}{2} \left( 1 + \frac{2 - 2 \delta}{1 + \lambda} - \sqrt{\left( 1 + \frac{2 - 2 \delta}{1 + \lambda} \right)^2 - \frac{8}{1 + \lambda}} \right),
\end{equation}
and $\omega_{i_0 + i} := \omega_{i_0} r(\delta)^i$ for $i \geq 1$ such that
\begin{equation}
    \omega_{i + 1} = \omega_i + \frac{2 \omega_i (1 - \delta) - 2 \omega_{i - 1}}{1 + \lambda} \leq \omega_i + \frac{\omega_i (1 - \delta) - \omega_{i - 1}}{\lambda (2 q_i^*(\infty) + 1)},
\end{equation}
for $i \geq i_0 + 1$. Note that $r(\delta) \to 1$ as $\delta \to 0$. Hence, there exists $\delta_2 > 0$ small such that $1 < r(\delta) < \min(2 / (1 + \lambda), 1 + \epsilon)$ for $\delta \leq \delta_2$. Finally, let $\delta := \min(\delta_0, \delta_1, \delta_2)$ and $r = r(\delta)$, such that equation \eqref{eq:globalstabcon} is satisfied.
\end{proof}
\section{Proportional sparsity and quasi-randomness}
\label{app:proofquasirandom}


\begin{proof}[Proof of Lemma~\ref{lem:quasi-to-prop}]
Note that the quasi-random graph can be characterized by the discrepancy condition \cite[Property 4]{CGW89}, which for bipartite graphs, states that a graph sequence $\{G_N\}_{N\geq 1}$ is quasi-random if there exists a fixed $0 < p < 1$, such that for all $U_N \subseteq V_N$ and $B_N \subseteq W_N$:
\begin{equation}
\label{eq:discrepancycon}
    \left\lvert \sum_{w \in B_N} \lvert \cN_w \cap U_N \rvert - p \lvert U_N \rvert \lvert B_N \rvert \right\rvert = o(N M(N)).
\end{equation}
We proceed by contradiction. Assume if possible, that there exists $\varepsilon > 0$ and a sequence of choices for $N$ for which there exists $U \subseteq V_N$ such that
\begin{equation}
    \left\lvert \left\{ w \in W_N \mid \left\lvert \frac{\lvert \cN_w \cap U \rvert}{\lvert \cN_w \rvert} - \frac{\lvert U \rvert}{N} \right\rvert \geq \varepsilon \right\} \right\rvert \geq 2 \delta M(N),
\end{equation}
for $\delta > 0$. Define the over connected set of dispatchers as
\begin{equation}
    B_1 := \left\{ w \in W_N \mid \frac{\lvert \cN_w \cap U \rvert}{\lvert \cN_w \rvert} - \frac{\lvert U \rvert}{N} \geq \varepsilon \right\},
\end{equation}
and assume without loss of generality that $\lvert B_1 \rvert \geq \delta M(N)$. Let $p := \sum_{w \in B_1} \lvert \cN_w \rvert / (\lvert B_1 \rvert N)$ be the average connection probability. Now,
\begin{equation}
    \sum_{w \in B_1} \lvert \cN_w \cup U \rvert \geq \left( \frac{\lvert U \rvert}{N} + \varepsilon \right) \sum_{w \in B_1} \left\lvert \cN_w \right\rvert
    = p \lvert U \rvert \lvert B_1 \rvert + \varepsilon p \lvert B_1 \rvert N
    \geq p \lvert U \rvert \lvert B_1 \rvert + \varepsilon \delta p N M(N).
\end{equation}
Define the under connected set of dispatchers to $U^c$ as
\begin{equation}
\begin{split}
    B_2 :&= \left\{ w \in W_N \Big| \frac{\lvert \cN_w \cap U^c \rvert}{\lvert \cN_w \rvert} - \frac{\lvert U^c \rvert}{N} \leq -\varepsilon \right\} \\
    &= \left\{ w \in W_N \Big| \left( 1 - \frac{\lvert \cN_w \cap U \rvert}{\lvert \cN_w \rvert} \right) - \left( 1 - \frac{\lvert U \rvert}{N} \right) \leq -\varepsilon \right\} \\
    &= \left\{ w \in W_N \Big| \frac{\lvert \cN_w \cap U \rvert}{\lvert \cN_w \rvert} - \frac{\lvert U \rvert}{N} \geq \varepsilon \right\} = B_1.
\end{split}
\end{equation}
With the same reasoning as before
\begin{equation}
    \sum_{w \in B_2} \lvert \cN_w \cup U^c \rvert \leq \left( \frac{\lvert U^c \rvert}{N} - \varepsilon \right) \sum_{w \in B_2} \left\lvert \cN_w \right\rvert
    \leq p \lvert U^c \rvert \lvert B_1 \rvert - \varepsilon \delta p N M(N).
\end{equation}
As the deviations are of the order $N M(N)$, it now follows that there is no choice for $p$ for which equation \eqref{eq:discrepancycon} is satisfied.
\end{proof}
\section{Lipschitz continuity of JSQ(\texorpdfstring{$d$}{d})}
\label{app:prooflemmajsqdlipschitz}


\begin{proof}[Proof of Lemma \ref{lem:jsqdlipschitz}]
Let $\xx \in \cX$ be the queue length distribution and $q_i = \sum_{j = i}^\infty x_j$ be the corresponding occupancy process. Let $\Pi$ be the JSQ($d$) policy. The assignment probability function $\pp^\Pi = \left( p_0^\Pi, p_1^\Pi, \dots \right) : \cX \to [0, 1]^\infty$ is given by
\begin{equation}
    p_{i-1}^\Pi(\xx) = q_{i-1}^d - q_i^d = (x_i + q_i)^d - q_i^d = \sum_{k = 0}^{d-1} \binom{d}{k} x_i^{d-k} q_i^k \quad\text{for } i = 1, 2, \dots.
\end{equation}
Let $\yy \in \cX$ be another queue length distribution with corresponding occupancy process $r_i = \sum_{j = i}^\infty y_j$. By applying the triangle inequality and Lipschitz continuity of $f(z) = z^k$,
\begin{align*}
    &\sum_{i = 1}^\infty \left\lvert p_{i-1}^\Pi(\yy) - p_{i-1}^\Pi(\xx) \right\rvert \leq \sum_{i = 1}^\infty \sum_{k = 0}^{d-1} \binom{d}{k} \left\lvert y_i^{d-k} r_i^k - x_i^{d-k} q_i^k \right\rvert \\
    &\leq d! \sum_{i = 1}^\infty \sum_{k = 0}^{d-1} \left( \left\lvert y_i^{d-k} r_i^k - x_i^{d-k} r_i^k \right\rvert + \left\lvert r_i^k - q_i^k \right\rvert \left\lvert x_i^{d-k} \right\rvert \right) \\
    &\leq d! \sum_{k = 0}^{d-1} \sum_{i = 1}^\infty \left\lvert y_i^{d-k} - x_i^{d-k} \right\rvert + d! \sum_{k = 0}^{d-1} \max_{i \geq 1} \left\lvert r_i^k - q_i^k \right\rvert \sum_{i = 1}^\infty \left\lvert x_i^{d-k} \right\rvert \\
    &\leq d! \sum_{k = 0}^{d-1} (d - k) \sum_{i = 1}^\infty \left\lvert y_i - x_i \right\rvert + d! \cdot d^2 \max_{i \geq 1} \left\lvert r_i - q_i \right\rvert \\
    &\leq d! \cdot d^2 \sum_{i = 0}^\infty \left\lvert y_i - x_i \right\rvert + d! \cdot d^2 \max_{i \geq 1} \left\lvert \sum_{j = i}^\infty \left( x_j - y_j \right) \right\rvert
    \leq 2 d! \cdot d^2 \sum_{i = 0}^\infty \left\lvert y_i - x_i \right\rvert.
\end{align*}

\end{proof}
\section{Proof of Proposition ~\ref{lemma:queuedistr}}
\label{app:queuedistr}

To prove Proposition~\ref{lemma:queuedistr}, we first need to show that the tail of the occupancy process is small uniformly on any finite time interval, for all large enough $N$. 
This is stated in the next lemma. 

\begin{lemma}
\label{lemma:queuetail}
If the starting states satisfy $\|\qq(\Phi(G_N, 0)) - \qq^*(0) \|_1 \to 0$ as $N \to \infty$, for some $\qq^*(0) \in \cY$, then for each $\varepsilon > 0$, $\delta > 0$, and $T > 0$ there exist $i_0 \geq 1$ and $N_1\geq 1$, possibly depending on $\lambda$, $\qq^*(0)$, $\varepsilon$, $\delta$, and $T$, such that,
\begin{equation}
    \mathbb{P}\left( \sup_{t \in [0, T]} q_{i_0}(\Phi(G_N, t)) \geq \varepsilon \right) < \delta \quad\text{for all}\quad N \geq N_1.
\end{equation}
\end{lemma}

Note that Lemma~\ref{lemma:queuetail} does not depend on any condition on the graph sequence nor does it require the Lipschitz continuity of the task assignment policy.

\begin{proof}[Proof of Lemma \ref{lemma:queuetail}]
Fix $\varepsilon > 0$ and $\delta > 0$. As $\qq^*(0) \in \ell_1$, there exists $j_0 = j_0(\qq^*(0)) \geq 1$ such that $q_{j_0}^*(0) < \varepsilon / 4$. By convergence of the starting states there exists $N_0 \geq 1$ such that
\begin{equation}
    \mathbb{P}(q_{j_0}(\Phi(G_N, 0)) \geq \varepsilon / 2) \leq \mathbb{P}(\left\lVert \qq(\Phi(G_N, 0)) - \qq^*(0) \right\rVert_1 \geq \varepsilon / 4) < \delta / 2 \quad\text{for}\quad N \geq N_0.
\end{equation}
Let $i_0 = j_0 + \lceil 4 \lambda T / \varepsilon \rceil$. Then,
\begin{equation}
\begin{gathered}
    \mathbb{P}\left( \sup_{t \in [0, T]} q_{i_0}(\Phi(G_N, t)) \geq \varepsilon \right) \\
    \leq \mathbb{P}\left( \sup_{t \in [0, T]} q_{i_0}(\Phi(G_N, t)) \geq \varepsilon \;\Big|\; q_{j_0}(\Phi(G_N, 0)) < \varepsilon/2 \right) + \mathbb{P}(q_{j_0}(\Phi(G_N, 0)) \geq \varepsilon/2).
\end{gathered}
\end{equation}
By the choice of $j_0$, the second term is smaller than $\delta / 2$ for $N \geq N_0$. By the condition $q_{j_0}(\Phi(G_N, 0)) < \varepsilon/2$, the number of servers with at least $j_0$ tasks is upper bounded by $\varepsilon N / 2$ at time zero. Hence to reach $q_{i_0}(\Phi(G_N, t)) \geq \varepsilon$, there must be at least $(\varepsilon N - \varepsilon N / 2) (i_0 - j_0) \geq 2 \lambda N T$ tasks added to the system. Therefore there exists $N_0' \geq 1$ such that
\begin{equation}
\begin{gathered}
    \mathbb{P}\left( \sup_{t \in [0, T]} q_{i_0}(\Phi(G_N, t)) \geq \varepsilon \;\Big|\; q_{j_0}(\Phi(G_N, 0)) < \varepsilon/2 \right)
    \leq \mathbb{P}(Z(\lambda N T) \geq 2 \lambda N T) < \delta / 2
\end{gathered}
\end{equation}
for $N \geq N_0'$,
where $Z(\cdot)$ is a unit-rate Poisson process and $Z(\lambda N T)$ denotes the total number of arrivals into the system up to time $T$. 
Choosing $N_1 = \max(N_0, N_0')$ completes the proof of the lemma.
\end{proof}

\begin{proof}[Proof of Proposition \ref{lemma:queuedistr}]
We can upper bound the probability in equation \eqref{eq:badsetsmall} by repeatedly applying the triangle inequality and splitting the probabilities in separate terms. For brevity, let $\xx(t) := \xx(\Phi(G_N, t))$ and $\xx^w(t) := \xx^w(\Phi(G_N, t))$ and similarly $\qq(t) := \qq(\Phi(G_N, t))$ and $\qq^w(t) := \qq^w(\Phi(G_N, t))$. Then,
\begin{equation}
\begin{split}
&\mathbb{P}\left( \sup_{t \in [0, T]} \left\lvert \left\{ w \in W_N \mid \sum_{i=0}^\infty \left\lvert x_i(t) - x_i^w(t) \right\rvert \geq \varepsilon \right\} \right\rvert \geq \delta M(N) \right) \\
&\leq \mathbb{P}\Bigg( \sup_{t \in [0, T]} \left\lvert \left\{ w \in W_N \mid \sum_{i=0}^{i_0 - 1} \left\lvert x_i(t) - x_i^w(t) \right\rvert \geq \varepsilon / 4 \right\} \right\rvert \\
&
\hspace{3cm}+ \sup_{t \in [0, T]} \left\lvert \left\{ w \in W_N \mid \sum_{i=i_0}^\infty x_{i_0}^w(t) \geq \varepsilon / 2 \right\} \right\rvert\\
&\hspace{4cm}+ \sup_{t \in [0, T]} \left\lvert \left\{ w \in W_N \mid \sum_{i=i_0}^\infty x_{i_0}(t) \geq \varepsilon / 4 \right\} \right\rvert
\geq \delta M(N) \Bigg) 
\end{split}
\end{equation}
\begin{equation}
\begin{split}\label{eq:prop2.6-1}
&
\leq \mathbb{P}\left( \sup_{t \in [0, T]} \left\lvert \left\{ w \in W_N \mid \sum_{i=0}^{i_0 - 1} \left\lvert x_i(t) - x_i^w(t) \right\rvert \geq \varepsilon / 4 \right\} \right\rvert > \delta M(N) / 4 \right) \\
&
\hspace{3cm}+ \mathbb{P}\left( \sup_{t \in [0, T]} \left\lvert \left\{ w \in W_N \mid q_{i_0}^w(t) \geq \varepsilon / 2 \right\} \right\rvert > \delta M(N) / 2 \right)\\
&\hspace{5cm}+ \mathbb{P}\left( \sup_{t \in [0, T]} \lvert \left\{ w \in W_N \mid q_{i_0}(t) \geq \varepsilon / 4 \right\} \rvert > \delta M(N) / 4 \right) \\
&
\leq \sum_{i=0}^{i_0 - 1} \mathbb{P}\left( \sup_{t \in [0, T]} \left\lvert \left\{ w \in W_N \mid \left\lvert x_i(t) - x_i^w(t) \right\rvert \geq \varepsilon / 4 i_0 \right\} \right\rvert > \delta M(N) / 4 i_0 \right) \\
&
\hspace{3cm}+ \mathbb{P}\left( \sup_{t \in [0, T]} \left\lvert \left\{ w \in W_N \mid \left\lvert q_{i_0}(t) - q_{i_0}^w(t) \right\rvert \geq \varepsilon / 4 \right\} \right\rvert > \delta M(N) / 4 \right)\\
&\hspace{5cm}
+ 2 \mathbb{P}\left( \sup_{t \in [0, T]} \mathbbm{1}\{ q_{i_0}(t) \geq \varepsilon / 4 \} > 0 \right),
\end{split}
\end{equation}
for an arbitrary choice of $i_0 \geq 1$, where in the second step we use $q_{i_0}(t) = \sum_{i = i_0}^\infty x_i(t)$. By Markov's inequality,
\begin{equation}
\begin{split}\label{eq:prop2.6-2b}
    \sum_{i = 0}^{i_0 - 1} \mathbb{P}\left( \sup_{t \in [0, T]} \left\lvert \left\{ w \in W_N \mid \left\lvert x_i(t) - x_i^w(t) \right\rvert \geq \varepsilon / 4 i_0 \right\} \right\rvert > \delta M(N) / 4 i_0 \right) \\
    \leq \frac{4 i_0}{\delta M(N)} \sum_{i = 0}^{i_0 - 1} \mathbb{E}\left[ \sup_{t \in [0, T]} \left\lvert \left\{ w \in W_N \mid \left\lvert x_i(t) - x_i^w(t) \right\rvert \geq \varepsilon / 4 i_0 \right\} \right\rvert \right] \\
    \leq \frac{4 i_0^2}{\delta M(N)} \sup_{U \subseteq V_N} \left\lvert \left\{ w \in W_N \mid \left\lvert \frac{\lvert \cN_w \cap U \rvert}{\lvert \cN_w \rvert} - \frac{\lvert U \rvert}{N} \right\rvert \geq \varepsilon / 4 i_0 \right\} \right\rvert,
\end{split}
\end{equation}
where the last step follows because the servers with queue length $i$ form a subset of $V_N$ at every time $t \geq 0$. Similarly, the second term in the right-hand side of \eqref{eq:prop2.6-1} can be bounded by
\begin{equation}
\begin{split}\label{eq:prop2.6-2}
    \mathbb{P}\left( \sup_{t \in [0, T]} \left\lvert \left\{ w \in W_N \mid \left\lvert q_{i_0}(t) - q_{i_0}^w(t) \right\rvert \geq \varepsilon / 4 \right\} \right\rvert > \delta M(N) / 4 \right) \\
    \leq \frac{4}{\delta M(N)} \mathbb{E}\left[ \sup_{t \in [0, T]} \left\lvert \left\{ w \in W_N \mid \left\lvert q_{i_0}(t) - q_{i_0}^w(t) \right\rvert \geq \varepsilon / 4 \right\} \right\rvert \right] \\
    \leq \frac{4}{\delta M(N)} \sup_{U \subseteq V_N} \left\lvert \left\{ w \in W_N \mid \left\lvert \frac{\lvert \cN_w \cap U \rvert}{\lvert \cN_w \rvert} - \frac{\lvert U \rvert}{N} \right\rvert \geq \varepsilon / 4 \right\} \right\rvert.
\end{split}
\end{equation}
Fix $\varepsilon' > 0$. By Lemma \ref{lemma:queuetail} we can choose $i_0 \geq 1$ and $N_0 \geq 1$ such that
\begin{equation}\label{eq:prop2.6-3}
    \mathbb{P}\left( \sup_{t \in [0, T]} q_{i_0}(t) \geq \varepsilon / 4 \right) \leq \varepsilon' / 4 \quad\text{for}\quad N \geq N_0.
\end{equation}
Since the graph sequence is proportionally sparse, by Definition~\ref{con:graphseq}, there exists $N_1 \geq 1$ such that
\begin{equation}
\begin{gathered}\label{eq:prop2.6-4}
    \sup_{U \subseteq V_N} \left\lvert \left\{ w \in W_N \mid \left\lvert \frac{\lvert \cN_w \cap U \rvert}{\lvert \cN_w \rvert} - \frac{\lvert U \rvert}{N} \right\rvert \geq \varepsilon / 4 i_0 \right\} \right\rvert < \frac{\delta}{4 i_0^2 + 4} \frac{\varepsilon' M(N)}{2} \quad\text{for all}\quad N \geq N_1.
\end{gathered}
\end{equation}
Therefore, plugging the bound from \eqref{eq:prop2.6-4} in \eqref{eq:prop2.6-2} and \eqref{eq:prop2.6-2b} and plugging this together with \eqref{eq:prop2.6-3} in \eqref{eq:prop2.6-1}, we obtain
\begin{equation}
\begin{gathered}
    \mathbb{P}\left( \sup_{t \in [0, T]} B_N^\varepsilon(t) \geq \delta M(N) \right) \leq \frac{4 i_0^2 + 4}{\delta M(N)} \sup_{U \subseteq V_N} \left\lvert \left\{ w \in W_N \mid \left\lvert \frac{\lvert \cN_w \cap U \rvert}{\lvert \cN_w \rvert} - \frac{\lvert U \rvert}{N} \right\rvert \geq \varepsilon / 4 i_0 \right\} \right\rvert \\
    + 2 \mathbb{P}\left( \sup_{t \in [0, T]} q_{i_0}(t) \geq \varepsilon / 4 \right) < \varepsilon' / 2 + \varepsilon' / 2 = \varepsilon' \quad\text{for all}\quad N \geq \max(N_0, N_1).
\end{gathered}
\end{equation}
\end{proof}
\section{Proof of Proposition \ref{prop:qlessdelta}}
\label{app:qlessdelta}

\begin{remark}\normalfont
The statement of Proposition \ref{prop:qlessdelta} is similar, in spirit, to \cite[Proposition~4]{MBLW16-3}. Our definition of mismatch in queue length counts the number of times a task is routed to servers with different queue lengths. In contrast, Mukherjee et al. \cite{MBLW16-3} order the servers by the queue lengths and their definition of \emph{differing in decision} counts the number of times a task is routed to servers with different order statistics, irrespective of whether they have the same queue length or not. Our notion of mismatch in queue length enables us to compare the occupancy processes of two different systems based on their LEQDs only, even when the individual queues are not coupled.
\end{remark}

\begin{proof}[Proof of Proposition \ref{prop:qlessdelta}]
We prove the inequality in equation~\eqref{eq:qlessdelta} by induction on the event times. Assume the inequality holds before the current time epoch $t$. We continue by distinguishing two cases, depending on whether $t$ is an arrival or a departure epoch. 
Recall that the $G_N$-system and $K_{N,M}$-system refer to the systems where the compatibility graph is $G_N$ and $K_{N,M}$, respectively.
Also, recall that, due to the optimal coupling, the servers in both systems are ordered by non-decreasing queue lengths and departures happen simultaneously at the $j$-th ordered server in the two systems, whenever they are non-empty for $j = 1, 2 \ldots, N$. \\

First, assume $t$ to be a departure epoch at the $j$-th ordered server. Our goal is to bound the increase of $\sum_{i=1}^\infty \left\lvert Q_i(\Phi(K_{N,M}, t)) - Q_i(\Phi(G_N, t)) \right\rvert$ by zero, because $\Delta_N(t)$ remains unchanged during the departure.
Let $I_1$ and $I_2$ denote the queue lengths just before time $t$ at the $j$-th ordered server for the $G_N$-system and $K_{N,M}$-system, respectively. After the departure, both $Q_{I_1}(\Phi(G_N, t))$ and $Q_{I_2}(\Phi(K_{N,M}, t))$ decrease by one, while the rest of the terms in the left-hand side of \eqref{eq:qlessdelta} remain unchanged. We will consider three possibilities:
\begin{enumerate}
    \item If $I_1 = I_2 = i_0$, then $Q_{i_0}$ decreases by one in both systems and hence, the sum in the left-hand side of~\eqref{eq:qlessdelta} remains unchanged. 
    \item Now assume $I_2 < I_1$. This implies that $Q_{I_1}(\Phi(G_N, t-)) > Q_{I_1}(\Phi(K_{N,M}, t-))$ and hence $\lvert Q_{I_1}(\Phi(K_{N, M}, t)) - Q_{I_1}(\Phi(G_N, t)) \rvert$ decreases by one during the departure. Since $\lvert Q_{I_2}(\Phi(K_{N, M}, t)) - Q_{I_2}(\Phi(G_N, t)) \rvert$ increases by at most one, the sum in the left-hand side of~\eqref{eq:qlessdelta} remains unchanged.
    \item The case $I_2 > I_1$ is similar to the $I_2 < I_1$ case.
\end{enumerate}
Next, assume $t$ to be an arrival epoch.
In this case, we further distinguish two cases:
\begin{enumerate}
    \item If the arriving tasks are routed to servers with unequal queue lengths in the two systems, then there is a mismatch in queue length at time $t$ and $\Delta_N(t)$ increases by one. Let $I_1$ and $I_2$ denote the queue lengths of the server the tasks are routed to for the $G_N$-system and $K_{N,M}$-system, respectively. During the arrival, $Q_{I_1 + 1}(\Phi(G_N, t))$ and $Q_{I_2 + 1}(\Phi(K_{N,M}, t))$ increase by one, while the rest of the terms in the left-hand side of \eqref{eq:qlessdelta} remain unchanged. Hence, $\sum_{i=1}^\infty \left\lvert Q_i(\Phi(K_{N,M}, t)) - Q_i(\Phi(G_N, t)) \right\rvert$ increases by at most two and, because $\Delta_N(t)$ increases by one, the right-hand side of (\ref{eq:qlessdelta}) increases by two as well.
    \item If the arriving tasks are routed to servers with equal queue lengths in the two systems, then there is no mismatch in queue length at time $t$ and $\Delta_N(t)$ remains unchanged. Therefore, our goal in this case is to bound the increase in the sum $\sum_{i=1}^\infty \left\lvert Q_i(\Phi(K_{N,M}, t)) - Q_i(\Phi(G_N, t)) \right\rvert$ by zero. 
    Let $I_1 = I_2 = i_0$ denote the queue length of the servers the tasks are routed to. During the arrival, $Q_{i_0 + 1}$ increases by one in both systems and hence, $\lvert Q_{i_0 + 1}(\Phi(K_{N, M}, t)) - Q_{i_0 + 1}(\Phi(G_N, t)) \rvert$ remains unchanged. As the rest of the terms in the left-hand side of \eqref{eq:qlessdelta} remain unchanged as well, the increase in the left-hand side of (\ref{eq:qlessdelta}) is bounded by zero.\\
\end{enumerate}
Therefore, in all the above cases, the inequality in~\eqref{eq:qlessdelta} is preserved at time epoch $t$.
This completes the proof of Proposition \ref{prop:qlessdelta}.
\end{proof}
\section{Proof of stability}
\label{app:proofstability}
As mentioned in the introduction, to prove stability and tightness of the steady state occupancy process, we use the Lyapunov function approach, as in~\cite{WMSY18, WMSY18a}, and establish moment bounds~\cite{Hajek82, MT93} to obtain uniform bounds on the tail of the stationary occupancy process.
To apply techniques from \cite{MT93} in discrete time, we consider the state of the system at \emph{event times} $t_0 = 0 < t_1 < \dots$, that is, for all $i\geq 1$, $t_i$ is either an arrival or a potential departure epoch. The Markov process can now be viewed as a uniformized Markov chain. We will later relate the behavior of this uniformized Markov chain to the behavior of the original process.

For any real-valued function $V: \cS_N \to \mathbb{R}$ defined on the state space, denote the expected increase $\Delta V(\cdot)$ as
\begin{equation}
    \Delta V(z) := \mathbb{E}[V(\Phi(G_N, t_1)) - V(\Phi(G_N, t_0)) \mid \Phi(G_N, t_0) = z].
\end{equation}
We investigate positive Harris recurrence of the uniformized chain by employing the next theorem.

\begin{theorem}[{\cite[Theorem 11.0.1]{MT93}}]
\label{th:posharrisrecurrence}
Suppose that $\Phi$ is a Markov chain. The Markov chain $\Phi$ is a positive Harris recurrent chain if and only if there exists some petite set $C \subseteq \cS_N$ and some non-negative function $V: \cS_N \to \mathbb{R}$, satisfying
\begin{equation}
    \Delta V(z) \leq -1 + b \mathbbm{1} \{ z \in C \} \quad\text{for all}\quad z \in \cS_N .
\end{equation}
\end{theorem}
The goal is to find an appropriate petite set $C$ and a suitable Lyapunov function $V$. 

\subsection{The Lyapunov function}
\label{sec:lyapunov}

Define a sequence of Lyapunov functions $V_k: \cS_N \to \mathbb{R}$ indexed by $k = 1, 2, \dots$ as
\begin{equation}
\label{eq:lyapunov}
V_k(z) := \sum_{i = k}^\infty \sum_{j = i}^\infty Q_j(z) 
= \frac{1}{2}\sum_{i = k}^\infty (i - k + 1)(i - k + 2) \cdot X_i(z).
\end{equation}
The sequence of Lyapunov functions, instead of a single Lyapunov function, will be necessary to bound the tail of the occupancy of the stationary state. This in turn will establish that the steady state of the occupancy process is tight in the appropriate space.

\begin{lemma}
\label{lemma:lyapunovdrift}
Consider a sequence $\{ G_N \}_{N \geq 1}$ of graphs which satisfies the subcriticality condition. For each $\varepsilon > 0$, there exists $N_0 \geq 1$ such that under the JSQ($d$) policy, for all $N \geq N_0$,
\begin{equation}
    \Delta V_k(z) \leq \frac{1}{\lambda + 1} \left( (1 + \varepsilon) \lambda q_{k - 1}(z) - (1 - (1 + \varepsilon) \lambda) \sum_{i = k}^\infty q_i(z) \right) \quad\text{for all } z \in \cS_N.
\end{equation}
\end{lemma}

Note that due to Lemma~\ref{lemma:lyapunovdrift}, if the tail of the occupancy is heavy compared to $q_{k - 1}(z)$, then the increase of the Lyapunov function is negative. In other words, if queues are long, a task is more likely to join shorter queues than to join one of these long queues. This is a consequence of the JSQ($d$) policy, which prefers shorter queues. The subcriticality condition on the graph is necessary here to ensure tasks are able to reach the shorter queues.

\begin{proof}[Proof of Lemma \ref{lemma:lyapunovdrift}]
By conveniently employing Poisson thinning as an argument, we can formulate an alternative description of the queuing system. Tasks in the system arrive as a Poisson process with rate $\lambda N$. At the arrival of a task, a dispatcher $w \in W_N$ is uniformly chosen at random. Next, a subset $U \subseteq \cN_w$ of size $\lvert U \rvert = d$ is selected uniformly at random over all subsets of $\cN_w$ of size $d$. The task is then routed to the server with the shortest queue in $U$. In this formulation the probability of a task being routed to the subset of servers with at least $i$ tasks, $\cQ_i$, is
\begin{equation}
    \frac{1}{M(N)} \sum_{w \in W_N} \binom{\lvert \cN_w \rvert}{d}^{-1} \sum_{\substack{U \subseteq \cN_w \\ \lvert U \rvert = \min(d, \lvert \cN_w \rvert)}} \left( \mathbbm{1}\{ U \subseteq \cQ_i \} \cdot 1 + (1 - \mathbbm{1}\{ U \subseteq \cQ_i \}) \cdot 0 \right).
\end{equation}
Because the graph sequence satisfies the subcriticality condition, we know there exist a probability distribution $\gamma_w^U(v)$ on $v \in U$ for each $w \in W_N$ and $U \subseteq V_N$ such that
\begin{equation}
    \max_{v \in V_N} \frac{N}{M(N)} \sum_{w \in W_N} \binom{\lvert \cN_w \rvert}{d}^{-1} \sum_{\substack{U \subseteq \cN_w \\ \lvert U \rvert = \min(d, \lvert \cN_w \rvert)}} \gamma_w^U(v) \leq 1 + \epsilon,
\end{equation}
for $N$ large enough. If $U \subseteq \cQ_i$ then
\begin{equation}
    \sum_{v \in \cQ_i} \gamma_w^U(v) = \sum_{v \in U} \gamma_w^U(v) = 1.
\end{equation}
and hence,
\begin{equation}
\begin{split}
    \frac{1}{M(N)} &\sum_{w \in W_N} \binom{\lvert \cN_w \rvert}{d}^{-1} \sum_{\substack{U \subseteq \cN_w \\ \lvert U \rvert = \min(d, \lvert \cN_w \rvert)}} \left( \mathbbm{1}\{ U \subseteq \cQ_i \} \cdot 1 + (1 - \mathbbm{1}\{ U \subseteq \cQ_i \}) \cdot 0 \right) \\
    &\leq \frac{1}{M(N)} \sum_{w \in W_N} \binom{\lvert \cN_w \rvert}{d}^{-1} \sum_{\substack{U \subseteq \cN_w \\ \lvert U \rvert = \min(d, \lvert \cN_w \rvert)}} \sum_{v \in \cQ_i} \gamma_w^U(v) \\
    &= \sum_{v \in \cQ_i} \frac{1}{M(N)} \sum_{w \in W_N} \binom{\lvert \cN_w \rvert}{d}^{-1} \sum_{\substack{U \subseteq \cN_w \\ \lvert U \rvert = \min(d, \lvert \cN_w \rvert)}} \gamma_w^U(v) \leq \sum_{v \in \cQ_i} \frac{1 + \varepsilon}{N} = (1 + \varepsilon) q_i,
\end{split}
\end{equation}
where the last inequality follows from the subcriticality condition for $N$ large enough. The increase of the Lyapunov function is then bounded as
\begin{equation}
\begin{split}
    \Delta V_k(z) &:= \mathbb{E}[V_k(\Phi(G_N, t_1)) - V_k(\Phi(G_N, t_0)) \mid \Phi_{t_0} = z] \\
    &= \sum_{i = k}^\infty \mathbb{E}\left[ \sum_{j = i}^\infty \left( Q_j(\Phi(G_N, t_1)) - Q_j(\Phi(G_N, t_0)) \right) \mid \Phi_{t_0} = z \right] \\
    &=  \sum_{i = k}^\infty \left( \frac{\lambda N}{\lambda N + N} \mathbb{P}\left( \text{task routed to } \cQ_{i - 1}(z) \right) - \frac{N}{\lambda N + N} \frac{Q_i(z)}{N} \right) \\
    &\leq \frac{1}{\lambda + 1} \sum_{i = k}^\infty \left( (1 + \varepsilon) \lambda q_{i - 1}(z) - q_i(z) \right)\\
    &= \frac{1}{\lambda + 1} \Big( (1 + \varepsilon) \lambda q_{k - 1}(z) - (1 - (1 + \varepsilon) \lambda) \sum_{i = k}^\infty q_i(z) \Big).
\end{split}
\end{equation}
\end{proof}

\subsection{Proof of positive Harris recurrence}
\label{sec:positiverecurrence}

We will gradually progress towards proving positive Harris recurrence using Theorem \ref{th:posharrisrecurrence}, starting with the definition of an appropriate petite set. The set $C_\Gamma \subseteq \cS_N$ defined as
\begin{equation}\label{eq:c-gamma}
    C_\Gamma := \left\{ z \in \cS_N \mid \sum_{i = 1}^\infty Q_i(z) \leq \Gamma \right\}.
\end{equation}
Recall the definition of a petite set.
\begin{definition}[Petite set]\textit{
A set $C \subseteq \cS_N$ is a \emph{petite set} if and only if there exists a distribution $a = \{ a(n) \}$ on $\mathbb{N}$ and a non-trivial measure $\mu$ on $\cS_N$ such that
\begin{equation}
    \sum_{n = 0}^\infty P^n(z, A) a(n) \geq \mu(A)
\end{equation}
for all $z \in C$ and $A \subseteq \cS_N$, where $P^n(z, A)$ is the $n$-step transition probability from $z$ to $A$.}
\end{definition}
Lemma~\ref{lemma:petiteset} below states that $C_\Gamma$ introduced in~\eqref{eq:c-gamma} is a petite set.
The proof is fairly straightforward since $C_\Gamma$ is finite, however, we provide it for the sake of completeness.
\begin{lemma}
\label{lemma:petiteset}
The set $C_\Gamma \subseteq \cS_N$ defined as in~\eqref{eq:c-gamma} is petite for any $\Gamma \geq 0$.
\end{lemma}

\begin{proof}
To establish petiteness, it is sufficient to show that there exists a non-trivial measure $\mu$ on $\cS_N$ such that
\begin{equation}
\label{eq:petitemeasure}
    \sum_{n = 1}^\infty e^{-n} P^n(z, A) \geq \mu(A),
\end{equation}
for all $z \in C_\Gamma$, $A \subseteq \cS_N$, where $P^n(z, A)$ is the $n$-step transition probability from $z$ to $A$. Define $\mu$ as
\begin{equation}
    \mu(a) := \min_{z \in C_\Gamma} \sum_{n = 1}^\infty e^{-n} P^n(z, a) \text{ for } a \in \cS_N,
\end{equation}
and let $\mu(A) = \sum_{a \in A} \mu(a)$ which is well-defined because the state space is countable. Trivially equation (\ref{eq:petitemeasure}) is satisfied. Moreover, because each state $a \in \cS_N$ is reachable and $\lvert C_\Gamma \rvert < \infty$, the minimum will be non-zero and the measure is non-trivial.
\end{proof}

The next lemma establishes that the Lyapunov function defined in the previous section satisfies the condition of Theorem \ref{th:stability}.

\begin{lemma}
\label{lemma:lyapunovdriftcon}
Consider a sequence $\{ G_N \}_{N \geq 1}$ of graphs which satisfies the subcriticality condition. There exists $N_0 \geq 1$ such that for the JSQ($d$) policy, the function $V_1$ as defined in equation (\ref{eq:lyapunov}) satisfies
\begin{equation}
\Delta V_1(z) \leq -1 + 2 \mathbbm{1}\{z \in C\} \quad\text{for all}\quad z \in \cS_N,
\end{equation}
for $C = C_\Gamma$ with $\Gamma = 3 \frac{1 + \lambda}{1 - \lambda} N$ and all $N \geq N_0$.
\end{lemma}

\begin{proof}
Choose $\varepsilon$ small such that $(1 + \varepsilon) \lambda \leq (\lambda + 1) / 2 < 1$. By Lemma \ref{lemma:lyapunovdrift}, we know that
\begin{equation}
\begin{gathered}
    \Delta V_1(z) \leq \frac{1}{\lambda + 1} \Big( (1 + \varepsilon) \lambda - (1 - (1 + \varepsilon) \lambda) \sum_{i = 1}^\infty q_i(z) \Big),
\end{gathered}
\end{equation}
for $N$ large enough. Then for $z \in \cS_N$,
\begin{equation}
    \Delta V_1(z) \leq \frac{(1 + \varepsilon) \lambda}{\lambda + 1} \leq \frac{1}{2} \frac{\lambda + 1}{\lambda + 1} < -1 + 2.
\end{equation}
Also, for $z \notin C_\Gamma$ with $\Gamma = 3 \frac{1 + \lambda}{1 - \lambda} N$,
\begin{equation}
\begin{gathered}
    \Delta V_1(z) \leq \frac{(1 + \varepsilon) \lambda}{1 + \lambda} - \frac{1 - (1 + \varepsilon) \lambda}{1 + \lambda} \sum_{i = 1}^\infty q_i(z) < \frac{1}{2} \frac{1 + \lambda}{1 + \lambda} - \frac{1}{2} \frac{1 - \lambda}{1 + \lambda} \frac{\Gamma}{N} \leq -1.
\end{gathered}
\end{equation}

\end{proof}

Lemmas \ref{lemma:petiteset} and \ref{lemma:lyapunovdriftcon} together with Theorem \ref{th:posharrisrecurrence} imply positive Harris recurrence of the uniformized Markov chain. This implies that also the continuous time chain is positive Harris recurrent and has a stationary distribution.

\subsection{Proof of moment bound}

Finally, we prove the moment bound in Lemma \ref{lemma:momentbound}.

\begin{proof}[Proof of Lemma \ref{lemma:momentbound}]
By the definition of the steady state,
\begin{equation}
     \mathbb{E}\left[ \Delta V_k(\Phi(G_N, \infty)) \right] = 0.
\end{equation}

Choose $\varepsilon$ small such that $(1 + \varepsilon) \lambda \leq (1 + \lambda) / 2 < 1$. Using Lemma \ref{lemma:lyapunovdrift} we conclude for $N$ large enough,
\begin{equation}
    \mathbb{E} \left[ (1 + \varepsilon) \lambda q_{k - 1}(z) - (1 - (1 + \varepsilon) \lambda) \sum_{i = k}^\infty q_i(z) \right] \geq 0,
\end{equation}
which can be rewritten to show the lemma.
\end{proof}
\section{Criteria for tightness in the \texorpdfstring{$\ell_1^\omega$}{weighted-l1}-topology}
\label{app:tightnessequiv}


\begin{proof}[Proof of Lemma \ref{lemma:tightnessequiv}]
Define the map $F : \ell_1 \rightarrow \ell_1^\omega$ as
\begin{equation}
    F(\boldsymbol{s}) = \left( \frac{s_1}{\omega_1}, \frac{s_2}{\omega_2}, \dots \right).
\end{equation}
$F$ is a bijective isometry as for $\boldsymbol{s}, \boldsymbol{t} \in \ell_1$,
\begin{equation}
    \lvert F(\boldsymbol{s}) - F(\boldsymbol{t}) \rvert_1^\omega = \sum_{i = 1}^\infty \omega_i \left\lvert \frac{s_i}{\omega_i} - \frac{t_i}{\omega_i} \right\rvert = \sum_{i = 1}^\infty \lvert s_i - t_i \rvert = \lVert \boldsymbol{s} - \boldsymbol{t} \rVert_1.
\end{equation}
For convenience, we denote $\boldsymbol{\omega} \cdot \boldsymbol{s} := F^{-1}(\boldsymbol{s})$. 
We first claim that the tightness of the sequence $\{ \qq(\Phi(G_N, \infty)) \}_{N \geq 1}$ in $\ell_1^\omega$ is equivalent to tightness of $\{ \boldsymbol{\omega} \cdot \qq(\Phi(G_N, \infty)) \}_{N \geq 1}$ in $\ell_1$.
Indeed, if $\{ \boldsymbol{\omega} \cdot \qq(\Phi(G_N, \infty)) \}_{N \geq 1}$ is tight in $\ell_1$, then we can find a compact set $K \subseteq \ell_1$ such that,
\begin{equation}
\label{eq:compactepsilon}
    \mathbb{P}( \qq(\Phi(G_N, \infty)) \notin F(K)) = \mathbb{P}( \boldsymbol{\omega} \qq(\Phi(G_N, \infty)) \notin K) < \epsilon.
\end{equation}
As the mapping preserves compactness due to continuity, $F(K)$ is a compact set in $\ell_1^\omega$ and hence $\{ \qq(\Phi(G_N, \infty)) \}_{N \geq 1}$ is tight in $\ell_1^\omega$. Conversely, if $\{ \qq(\Phi(G_N, \infty)) \}_{N \geq 1}$ is tight in $\ell_1^\omega$, then there is a compact set $K \subseteq \ell_1^\omega$ such that equation (\ref{eq:compactepsilon}) holds and $F^{-1}(K)$ is compact in $\ell_1$. Hence, the claim is proved.

Now, by Lemma 2 of \cite{MBLW16-3}, tightness in $\ell_1$ is equivalent to showing tightness with respect to the product topology and~\eqref{eq:tightvanishingtail}.
This completes the proof of Lemma \ref{lemma:tightnessequiv}.
\end{proof}

\end{document}